\documentclass{article}
\usepackage{float, amsmath, amssymb, amsfonts, enumerate, vmargin, wasysym, accents, multicol, tikz, spverbatim, cjhebrew, amsthm, relsize, framed, epstopdf, caption, chngcntr, url, xcolor}
\usepackage[bookmarks=false]{hyperref}

\usetikzlibrary{decorations.markings,calc,arrows}
\usetikzlibrary{calc}

\hypersetup{colorlinks=true, linkcolor=violet, filecolor=blue, citecolor=violet, urlcolor=darkgray, citebordercolor=lightgray, filebordercolor=lightgray, linkbordercolor=lightgray}

\counterwithout{paragraph}{subsubsection} 
\counterwithin{paragraph}{section}

\setmarginsrb{30mm}{20mm}{30mm}{20mm}{0pt}{10mm}{0pt}{10mm}

\def\N{\mathbb{N}} \def\Q{\mathbb{Q}} \def\Z{\mathbb{Z}} \def\R{\mathbb{R}}  \def\C{\mathbb{C}}

 \newtheorem*{thm*}{Theorem}
 \newtheorem*{claim*}{Claim}
 \newtheorem*{dfn*}{Definition}
 \newtheorem*{lemma*}{Lemma}
 \newtheorem*{prop*}{Proposition}
 \newtheorem*{cor*}{Corollary}
 \newtheorem*{conj*}{Conjecture}
 \newtheorem*{quest*}{Question}
 \newtheorem*{exer*}{Exercise}
\theoremstyle{remark}
 \newtheorem*{rmk*}{Remark}
 \newtheorem*{example*}{Example}
 \newtheorem*{examples*}{Examples}
 \newtheorem*{rmks*}{Remarks}

\newcommand{\knottable}{
\begin{tabular}{cccc} &&& \\
\begin{tikzpicture}
\draw[white,line width=2.0pt,double=black,double distance=2.0pt] (1,0) .. controls +(90:1.35) and +(90:1.35) .. (-1,0);
\draw[white,line width=2.0pt,double=black,double distance=2.0pt] (1,0) .. controls +(270:1.35) and +(270:1.35) .. (-1,0);
\draw[black,line width=2.0pt,decoration={markings,mark=at position 1 with {\arrow{>}}},postaction={decorate}] (1,0) -- +(90:0.1);
\pgfresetboundingbox \clip (-1.25,-1.25) rectangle (1.25,1.25);
\end{tikzpicture} &
\begin{tikzpicture}
\draw[white,line width=2.0pt,double=black,double distance=2.0pt] (90:0.4) .. controls +(180:0.6) and +(120:0.5) .. (210:1.2);
\draw[white,line width=2.0pt,double=black,double distance=2.0pt] (210:0.4) .. controls +(300:0.6) and +(240:0.5) .. (330:1.2);
\draw[white,line width=2.0pt,double=black,double distance=2.0pt] (330:0.4) .. controls +(60:0.6) and +(0:0.5) .. (90:1.2);
\draw[white,line width=2.0pt,double=black,double distance=2.0pt] (90:0.4) .. controls +(0:0.6) and +(60:0.5) .. (330:1.2);
\draw[white,line width=2.0pt,double=black,double distance=2.0pt] (210:0.4) .. controls +(120:0.6) and +(180:0.5) .. (90:1.2);

\draw[white,line width=2.0pt,double=black,double distance=2.0pt] (330:0.4) .. controls +(240:0.6) and +(300:0.5) .. (210:1.2);
\draw[black,line width=2.0pt,decoration={markings,mark=at position 1 with {\arrow{>}}},postaction={decorate}] (90:0.4) -- +(180:0.1);
\pgfresetboundingbox \clip (-1.25,-1.25) rectangle (1.25,1.25);
\end{tikzpicture} &
\begin{tikzpicture}
\draw[white,line width=2.0pt,double=black,double distance=2.0pt] (90:0.4) .. controls +(0:0.6) and +(60:0.5) .. (330:1.2);
\draw[white,line width=2.0pt,double=black,double distance=2.0pt] (210:0.4) .. controls +(120:0.6) and +(180:0.5) .. (90:1.2);
\draw[white,line width=2.0pt,double=black,double distance=2.0pt] (330:0.4) .. controls +(240:0.6) and +(300:0.5) .. (210:1.2);
\draw[white,line width=2.0pt,double=black,double distance=2.0pt] (90:0.4) .. controls +(180:0.6) and +(120:0.5) .. (210:1.2);
\draw[white,line width=2.0pt,double=black,double distance=2.0pt] (210:0.4) .. controls +(300:0.6) and +(240:0.5) .. (330:1.2);
\draw[white,line width=2.0pt,double=black,double distance=2.0pt] (330:0.4) .. controls +(60:0.6) and +(0:0.5) .. (90:1.2);
\draw[black,line width=2.0pt,decoration={markings,mark=at position 1 with {\arrow{>}}},postaction={decorate}] (90:0.4) -- +(0:0.1);
\pgfresetboundingbox \clip (-1.25,-1.25) rectangle (1.25,1.25);
\end{tikzpicture} &
\begin{tikzpicture}
\draw[white,line width=2.0pt,double=black,double distance=2.0pt] (0:0.5) .. controls +(90:0.4) and +(45:1.3) .. (135:1.2);
\draw[white,line width=2.0pt,double=black,double distance=2.0pt] (270:0.5) .. controls +(0:0.4) and +(315:0.2) .. (45:0.1);
\draw[white,line width=2.0pt,double=black,double distance=2.0pt] (180:0.5) .. controls +(270:0.4) and +(225:1.3) .. (315:1.2);
\draw[white,line width=2.0pt,double=black,double distance=2.0pt] (90:0.5) .. controls +(180:0.4) and +(135:0.2) .. (225:0.1);
\draw[white,line width=2.0pt,double=black,double distance=2.0pt] (135:1.2) .. controls +(225:1.3) and +(180:0.4) .. (270:0.5);
\draw[white,line width=2.0pt,double=black,double distance=2.0pt] (45:0.1) .. controls +(135:0.2) and +(90:0.4) .. (180:0.5);
\draw[white,line width=2.0pt,double=black,double distance=2.0pt] (315:1.2) .. controls +(45:1.3) and +(0:0.4) .. (90:0.5);
\draw[white,line width=2.0pt,double=black,double distance=2.0pt] (225:0.1) .. controls +(315:0.2) and +(270:0.4) .. (0:0.5);
\draw[black,line width=2.0pt,decoration={markings,mark=at position 1 with {\arrow{>}}},postaction={decorate}] (315:1.2) -- +(45:0.1);
\pgfresetboundingbox \clip (-1.25,-1.25) rectangle (1.25,1.25);
\end{tikzpicture} \\
Unknot & Left Trefoil & Right Trefoil & Figure Eight \\ &&& \\
\begin{tikzpicture}
\draw[white,line width=2.0pt,double=black,double distance=2.0pt] (1.1,0.5) .. controls +(90:0.8) and +(90:0.8) .. (-0.1,0.5);
\draw[white,line width=2.0pt,double=black,double distance=2.0pt] (1.1,0.5) .. controls +(270:0.8) and +(270:0.8) .. (-0.1,0.5);
\draw[white,line width=2.0pt,double=black,double distance=2.0pt] (0.1,-0.5) .. controls +(90:0.8) and +(90:0.8) .. (-1.1,-0.5);
\draw[white,line width=2.0pt,double=black,double distance=2.0pt] (0.1,-0.5) .. controls +(270:0.8) and +(270:0.8) .. (-1.1,-0.5);
\draw[black,line width=2.0pt,decoration={markings,mark=at position 1 with {\arrow{>}}},postaction={decorate}] (1.1,0.5) -- +(90:0.1);
\draw[black,line width=2.0pt,decoration={markings,mark=at position 1 with {\arrow{>}}},postaction={decorate}] (0.1,-0.5) -- +(90:0.1);
\pgfresetboundingbox \clip (-1.25,-1.25) rectangle (1.25,1.25);
\end{tikzpicture} &
\begin{tikzpicture}
\draw[white,line width=2.0pt,double=black,double distance=2.0pt] (0.3,0.3) .. controls +(315:1.0) and +(315:1.0) .. (-0.8,-0.8);
\draw[white,line width=2.0pt,double=black,double distance=2.0pt] (0.8,0.8) .. controls +(135:1.0) and +(135:1.0) .. (-0.3,-0.3);
\draw[white,line width=2.0pt,double=black,double distance=2.0pt] (0.3,0.3) .. controls +(135:1.0) and +(135:1.0) .. (-0.8,-0.8);
\draw[white,line width=2.0pt,double=black,double distance=2.0pt] (0.8,0.8) .. controls +(315:1.0) and +(315:1.0) .. (-0.3,-0.3);
\draw[black,line width=2.0pt,decoration={markings,mark=at position 1 with {\arrow{>}}},postaction={decorate}] (0.3,0.3) -- +(135:0.1);
\draw[black,line width=2.0pt,decoration={markings,mark=at position 1 with {\arrow{>}}},postaction={decorate}] (-0.3,-0.3) -- +(315:0.1);
\pgfresetboundingbox \clip (-1.25,-1.25) rectangle (1.25,1.25);
\end{tikzpicture} &
\begin{tikzpicture}
\draw[white,line width=2.0pt,double=black,double distance=2.0pt] (0.6,0) .. controls +(270:0.8) and +(0:0.5) .. (0,-1);
\draw[white,line width=2.0pt,double=black,double distance=2.0pt] (-0.6,0) .. controls +(90:0.8) and +(180:0.5) .. (0,1);
\draw[white,line width=2.0pt,double=black,double distance=2.0pt] (0,0) .. controls +(-45:0.4) and +(-90:0.8) .. (1.2,0);
\draw[white,line width=2.0pt,double=black,double distance=2.0pt] (0,0) .. controls +(135:0.4) and +(90:0.8) .. (-1.2,0);
\draw[white,line width=2.0pt,double=black,double distance=2.0pt] (0,0) .. controls +(45:0.4) and +(90:0.8) .. (1.2,0);
\draw[white,line width=2.0pt,double=black,double distance=2.0pt] (0,0) .. controls +(-135:0.4) and +(-90:0.8) .. (-1.2,0);
\draw[white,line width=2.0pt,double=black,double distance=2.0pt] (0.6,0) .. controls +(90:0.8) and +(0:0.5) .. (0,1);
\draw[white,line width=2.0pt,double=black,double distance=2.0pt] (-0.6,0) .. controls +(270:0.8) and +(180:0.5) .. (0,-1);
\draw[black,line width=2.0pt,decoration={markings,mark=at position 1 with {\arrow{>}}},postaction={decorate}] (1.2,0) -- +(90:0.1);
\draw[black,line width=2.0pt,decoration={markings,mark=at position 1 with {\arrow{>}}},postaction={decorate}] (0,1) -- +(180:0.05);
\pgfresetboundingbox \clip (-1.25,-1.25) rectangle (1.25,1.25);
\end{tikzpicture} &
\begin{tikzpicture}
\draw[white,line width=2.0pt,double=black,double distance=2.0pt] (-90:1.2) .. controls +(-0:0.6) and +(-60:0.6) .. (-0:0.7);
\draw[white,line width=2.0pt,double=black,double distance=2.0pt] (-210:1.2) .. controls +(-120:0.6) and +(-180:0.6) .. (-120:0.7);
\draw[white,line width=2.0pt,double=black,double distance=2.0pt] (-330:1.2) .. controls +(-240:0.6) and +(60:0.6) .. (-240:0.7);
\draw[white,line width=2.0pt,double=black,double distance=2.0pt] (-90:1.2) .. controls +(-180:0.6) and +(-120:0.6) .. (-180:0.7);
\draw[white,line width=2.0pt,double=black,double distance=2.0pt] (-210:1.2) .. controls +(-300:0.6) and +(120:0.6) .. (-300:0.7);
\draw[white,line width=2.0pt,double=black,double distance=2.0pt] (-330:1.2) .. controls +(-60:0.6) and +(0:0.6) .. (-60:0.7);
\draw[white,line width=2.0pt,double=black,double distance=2.0pt] (-90:-0.4) .. controls +(-180:0.3) and +(60:0.3) .. (-180:0.7);
\draw[white,line width=2.0pt,double=black,double distance=2.0pt] (-210:-0.4) .. controls +(-300:0.3) and +(-60:0.3) .. (-300:0.7);
\draw[white,line width=2.0pt,double=black,double distance=2.0pt] (-330:-0.4) .. controls +(-60:0.3) and +(180:0.3) .. (-60:0.7);
\draw[white,line width=2.0pt,double=black,double distance=2.0pt] (-90:-0.4) .. controls +(-0:0.3) and +(120:0.3) .. (-0:0.7);
\draw[white,line width=2.0pt,double=black,double distance=2.0pt] (-210:-0.4) .. controls +(-120:0.3) and +(0:0.3) .. (-120:0.7);
\draw[white,line width=2.0pt,double=black,double distance=2.0pt] (-330:-0.4) .. controls +(-240:0.3) and +(-120:0.3) .. (-240:0.7);
\draw[black,line width=2.0pt,decoration={markings,mark=at position 1 with {\arrow{>}}},postaction={decorate}] (-90:1.2) -- +(0:0.05);
\draw[black,line width=2.0pt,decoration={markings,mark=at position 1 with {\arrow{>}}},postaction={decorate}] (-210:1.2) -- +(-120:0.1);
\draw[black,line width=2.0pt,decoration={markings,mark=at position 1 with {\arrow{>}}},postaction={decorate}] (-330:1.2) -- +(-240:0.11);
\pgfresetboundingbox \clip (-1.25,-1.25) rectangle (1.25,1.25);
\end{tikzpicture} \\
Unlink & Hopf Link & Whitehead Link & Borromean Rings 
\end{tabular}}

\newcommand{\knotwalk}{
\begin{tabular}{cc} 
\begin{tikzpicture}
\foreach \i/\j/\k in {0/1/0,1/1/0,2/1/0,2/1/1,2/1/2,3/1/2,3/2/2,3/3/2,3/3/1,2/3/1,2/3/2,1/3/2,1/3/1,1/2/1,1/1/1,1/0/1,2/0/1,3/0/1,3/1/1,3/2/1,2/2/1,2/2/2,1/2/2,1/1/2,0/1/2,0/2/2,0/2/1,0/1/1,0/0/1,0/0/0} 
\node[fill, shading=ball, circle, ball color=gray] (\i\j\k) at (1.2*\i+0.5*\k,1.2*\j+0.3*\k) {};
\draw[white,line width=1.5pt,double=black,double distance=2.0pt] (221)--(222)--(122)--(112)--(012)--(022)--(021)--(011)--(001)--(000)--(010);
\draw[white,line width=1.5pt,double=black,double distance=2.0pt] (211)--(212)--(312)--(322)--(332)--(331)--(231)--(232)--(132)--(131)--(121);
\draw[white,line width=1.5pt,double=black,double distance=2.0pt] (121)--(111)--(101)--(201)--(301)--(311)--(321)--(221);
\draw[white,line width=1.5pt,double=black,double distance=2.0pt] (010)--(110)--(210)--(211);
\pgfresetboundingbox \clip (-1,-0.5) rectangle (6,5);
\end{tikzpicture} &
\begin{tikzpicture}
\foreach \i/\j/\k in {0/1/0,0/0/0,1/0/0,2/0/0,2/1/0,2/1/1,2/1/2,2/2/2,2/2/3,1/2/3,0/2/3,0/2/2,0/2/1,1/2/1,1/1/1,1/0/1,2/0/1,3/0/1,3/1/1,3/2/1,2/2/1,1/3/2,1/2/2,1/1/2,0/1/2,0/1/1,1/1/0,1/2/0,1/3/0,1/3/1} 
\node[fill, shading=ball, circle, ball color=gray] (\i\j\k) at (1.2*\i+0.5*\k,1.2*\j+0.3*\k) {};
\draw[white,line width=1.5pt,double=black,double distance=2.0pt] (212)--(222)--(223)--(123)--(023)--(022)--(021)--(121);
\draw[white,line width=1.5pt,double=black,double distance=2.0pt] (132)--(122)--(112)--(012)--(011)--(010);
\draw[white,line width=1.5pt,double=black,double distance=2.0pt] (110)--(111)--(101)--(201)--(301)--(311)--(321)--(221)--(121);
\draw[white,line width=1.5pt,double=black,double distance=2.0pt] (010)--(000)--(100)--(200)--(210)--(211)--(212);
\draw[white,line width=1.5pt,double=black,double distance=2.0pt] (110)--(120)--(130)--(131)--(132);
\pgfresetboundingbox \clip (-1,-0.5) rectangle (5,5);
\end{tikzpicture}
\end{tabular}}

\newcommand{\knotpoly}{
\begin{tikzpicture}
\foreach \i/\x/\y in {0/0.00/0.00,1/0.18/0.93,2/0.75/0.11,3/0.21/-0.20,4/0.26/0.42,5/1.41/0.52,6/1.84/-0.55,7/2.25/-0.21,8/1.35/-0.48,9/1.39/-0.02,10/1.55/-0.13,11/0.83/-0.56,12/0.94/0.24,13/1.15/0.86,14/0.38/0.13,15/0.00/0.00} 
\node[fill, shading=ball, circle, ball color=gray] (\i) at (3*\x,3*\y) {};
\draw[white,line width=1.5pt,double=black,double distance=2.0pt] (3)--(4);
\draw[white,line width=1.5pt,double=black,double distance=2.0pt] (7)--(8);
\draw[white,line width=1.5pt,double=black,double distance=2.0pt] (5)--(6);
\draw[white,line width=1.5pt,double=black,double distance=2.0pt] (6)--(7);
\draw[white,line width=1.5pt,double=black,double distance=2.0pt] (13)--(14);
\draw[white,line width=1.5pt,double=black,double distance=2.0pt] (14)--(15);
\draw[white,line width=1.5pt,double=black,double distance=2.0pt] (1)--(2);
\draw[white,line width=1.5pt,double=black,double distance=2.0pt] (0)--(1);
\draw[white,line width=1.5pt,double=black,double distance=2.0pt] (2)--(3);
\draw[white,line width=1.5pt,double=black,double distance=2.0pt] (10)--(11);
\draw[white,line width=1.5pt,double=black,double distance=2.0pt] (9)--(10);
\draw[white,line width=1.5pt,double=black,double distance=2.0pt] (8)--(9);
\draw[white,line width=1.5pt,double=black,double distance=2.0pt] (11)--(12);
\draw[white,line width=1.5pt,double=black,double distance=2.0pt] (12)--(13);
\draw[white,line width=1.5pt,double=black,double distance=2.0pt] (4)--(5);
\pgfresetboundingbox \clip (-1,-2) rectangle (7,3.5);
\end{tikzpicture}}

\newcommand{\RMIa}{
\tikz{
\draw[white,line width=1.2pt,double=black,double distance=1.2pt] (-0.25,0) .. controls +(0:0.1) and +(180:0.1) .. (0,0.05);
\draw[white,line width=1.2pt,double=black,double distance=1.2pt] (0,0.05) .. controls +(0:0.1) and +(180:0.1) .. (0.25,0);
\pgfresetboundingbox \clip (-0.25,0) rectangle (0.25,0.3);}}

\newcommand{\RMIb}{
\tikz{
\draw[white,line width=1.2pt,double=black,double distance=1.2pt] (-0.25,0) .. controls +(0:0.4) and +(0:0.15) .. (0,0.25);
\draw[white,line width=1.2pt,double=black,double distance=1.2pt] (0,0.25) .. controls +(180:0.15) and +(180:0.4) .. (0.25,0);
\pgfresetboundingbox \clip (-0.25,0) rectangle (0.25,0.3);}}

\newcommand{\RMIIa}{
\tikz{
\draw[white,line width=1.2pt,double=black,double distance=1.2pt] (-0.25,0) .. controls +(60:0.1) and +(120:0.1) .. (0.25,0);
\draw[white,line width=1.2pt,double=black,double distance=1.2pt] (-0.25,0.25) .. controls +(-60:0.1) and +(-120:0.1) .. (0.25,0.25);
\pgfresetboundingbox \clip (-0.25,0) rectangle (0.25,0.3);}}

\newcommand{\RMIIb}{
\tikz{
\draw[white,line width=1.2pt,double=black,double distance=1.2pt] (-0.25,0) .. controls +(60:0.3) and +(120:0.3) .. (0.25,0);
\draw[white,line width=1.2pt,double=black,double distance=1.2pt] (-0.25,0.25) .. controls +(-60:0.3) and +(-120:0.3) .. (0.25,0.25);
\pgfresetboundingbox \clip (-0.25,0) rectangle (0.25,0.3);}}

\newcommand{\RMIIIa}{
\tikz{
\draw[white,line width=1.2pt,double=black,double distance=1.2pt] (-0.3,0.125) .. controls +(0:0.1) and +(180:0.1) .. (0,0.25);
\draw[white,line width=1.2pt,double=black,double distance=1.2pt] (0,0.25) .. controls +(0:0.1) and +(180:0.1) .. (0.3,0.125);
\draw[white,line width=1.2pt,double=black,double distance=1.2pt] (-0.2,0.25) .. controls +(-90:0) and +(90:0) .. (0.15,-0.05);
\draw[white,line width=1.2pt,double=black,double distance=1.2pt] (-0.15,-0.05) .. controls +(90:0) and +(-90:0) .. (0.2,0.25);
\pgfresetboundingbox \clip (-0.3,0) rectangle (0.3,0.3);}}

\newcommand{\RMIIIb}{
\tikz{
\draw[white,line width=1.2pt,double=black,double distance=1.2pt] (-0.3,0.125) .. controls +(0:0.1) and +(180:0.1) .. (0,0);
\draw[white,line width=1.2pt,double=black,double distance=1.2pt] (0,0) .. controls +(0:0.1) and +(180:0.1) .. (0.3,0.125);
\draw[white,line width=1.2pt,double=black,double distance=1.2pt] (-0.15,0.3) .. controls +(-90:0) and +(90:0) .. (0.2,0);
\draw[white,line width=1.2pt,double=black,double distance=1.2pt] (-0.2,0) .. controls +(90:0) and +(-90:0) .. (0.15,0.3);
\pgfresetboundingbox \clip (-0.3,0) rectangle (0.3,0.3);}}

\newcommand{\trefoil}{
\tikz{
\draw[white,line width=0.5pt,double=black,double distance=0.5pt] (90:0.08) .. controls +(180:0.12) and +(120:0.1) .. (210:0.24);
\draw[white,line width=0.5pt,double=black,double distance=0.5pt] (210:0.08) .. controls +(300:0.12) and +(240:0.1) .. (330:0.24);
\draw[white,line width=0.5pt,double=black,double distance=0.5pt] (330:0.08) .. controls +(60:0.12) and +(0:0.1) .. (90:0.24);
\draw[white,line width=0.5pt,double=black,double distance=0.5pt] (90:0.08) .. controls +(0:0.12) and +(60:0.1) .. (330:0.24);
\draw[white,line width=0.5pt,double=black,double distance=0.5pt] (210:0.08) .. controls +(120:0.12) and +(180:0.1) .. (90:0.24);
\draw[white,line width=0.5pt,double=black,double distance=0.5pt] (330:0.08) .. controls +(240:0.12) and +(300:0.1) .. (210:0.24);
\draw[black,line width=0.5pt,decoration={markings,mark=at position 1 with {\arrow{>}}},postaction={decorate}] (90:0.08) -- +(180:0.02);
\pgfresetboundingbox \clip (-0.3,-0.05) rectangle (0.3,0.25);}}

\newcommand{\granny}{
\tikz{
\draw[white,line width=0.5pt,double=black,double distance=0.5pt] (90:0.08) .. controls +(180:0.12) and +(120:0.1) .. (210:0.24);
\draw[white,line width=0.5pt,double=black,double distance=0.5pt] (210:0.08) .. controls +(300:0.12) and +(240:0.1) .. (330:0.24);
\draw[white,line width=0.5pt,double=black,double distance=0.5pt] (330:0.08) .. controls +(60:0.12) and +(0:0.1) .. (90:0.24);
\draw[white,line width=0.5pt,double=black,double distance=0.5pt] (90:0.08) .. controls +(0:0.12) and +(90:0.04) .. (340:0.2); 
\draw[white,line width=0.5pt,double=black,double distance=0.5pt] (210:0.08) .. controls +(120:0.12) and +(180:0.1) .. (90:0.24);
\draw[white,line width=0.5pt,double=black,double distance=0.5pt] (330:0.08) .. controls +(240:0.12) and +(300:0.1) .. (210:0.24);
\draw[black,line width=0.5pt,decoration={markings,mark=at position 1 with {\arrow{>}}},postaction={decorate}] (90:0.08) -- +(180:0.02);
\begin{scope}[shift={(0.7,0)}]
\draw[white,line width=0.5pt,double=black,double distance=0.5pt] (90:0.08) .. controls +(180:0.12) and +(90:0.04) .. (200:0.2); 
\draw[white,line width=0.5pt,double=black,double distance=0.5pt] (210:0.08) .. controls +(300:0.12) and +(240:0.1) .. (330:0.24);
\draw[white,line width=0.5pt,double=black,double distance=0.5pt] (330:0.08) .. controls +(60:0.12) and +(0:0.1) .. (90:0.24);
\draw[white,line width=0.5pt,double=black,double distance=0.5pt] (90:0.08) .. controls +(0:0.12) and +(60:0.1) .. (330:0.24);
\draw[white,line width=0.5pt,double=black,double distance=0.5pt] (210:0.08) .. controls +(120:0.12) and +(180:0.1) .. (90:0.24); 
\draw[white,line width=0.5pt,double=black,double distance=0.5pt] (330:0.08) .. controls +(240:0.12) and +(300:0.1) .. (210:0.24);
\draw[black,line width=0.5pt,decoration={markings,mark=at position 1 with {\arrow{>}}},postaction={decorate}] (90:0.08) -- +(180:0.02);
\end{scope}
\draw[white,line width=0.5pt,double=black,double distance=0.5pt] (0.2,-0.06) .. controls +(30:0.1) and +(150:0.1) .. (0.5,-0.06);
\draw[white,line width=0.5pt,double=black,double distance=0.5pt] (0.21,-0.12) .. controls +(30:0.1) and +(150:0.1) .. (0.49,-0.12);
\pgfresetboundingbox \clip (-0.3,-0.05) rectangle (1,0.25);}}

\newcommand{\knotpetal}{
\tikz[thick,decoration={markings,mark=at position 0.05 with {\arrow{<}}}]{
\foreach \angle in {0, 40, ..., 320} \draw[postaction={decorate}] (\angle:1) .. controls +(\angle:1.5) and +(\angle+20:1.5) .. (\angle+20:1);
\foreach \angle in {0, 40, ..., 320} \draw (\angle:1) -- (\angle:-1);
\foreach \angle/\num in {0/0, 80/1, 160/2, 120/3, 240/4, 280/5, 200/6, 320/7, 40/8} \node[scale=0.7] at (-\angle+6:0.8) {$ $}; 
\pgfresetboundingbox \clip (-2.25,-2.25) rectangle (2.25,2.25);}}

\newcommand{\linkpetal}{
\tikz[thick,decoration={markings,mark=at position 0.5 with {\arrow{>}}},scale=0.925]{
\foreach \angle/\color in {30/black,50/black,120/gray,140/gray,200/black,220/black,240/black,290/gray,310/gray,330/gray} \draw[color=\color,line width=1] (0,0) .. controls +(\angle:3) and +(\angle+10:3) .. (0,0);
\foreach \angle/\color in {15/black,105/gray} \draw[postaction={decorate},color=\color,line width=1] (0,0) .. controls +(\angle:4.4) and +(\angle+60:4.4) .. (0,0);
\pgfresetboundingbox \clip (-2.25,-2.25) rectangle (2.25,2.25);}}

\newcommand{\trilinkpetal}{
\begin{tikzpicture}[thick,decoration={markings,mark=at position 0.5 with {\arrow{>}}}]
\foreach \angle/\color in {
202/black,212/black,222/black,192/black,
282/lightgray,262/lightgray,252/lightgray,272/lightgray,
342/gray,322/gray,312/gray,332/gray} 
\draw[color=\color] (0,0) .. controls +(\angle:2.4) and +(\angle+6:2.4) .. (0,0);
\foreach \angle/\color in {
17/black,27/black,37/black,
97/lightgray,87/lightgray,77/lightgray,
157/gray,147/gray,137/gray} 
\draw[color=\color] (0,0) .. controls +(\angle:2.2) and +(\angle+6:2.2) .. (0,0);
\foreach \angle/\color in {7/black,67/lightgray,127/gray} \draw[postaction={decorate},color=\color] (0,0) .. controls +(\angle:3.38) and +(\angle+46:3.38) .. (0,0);
\pgfresetboundingbox \clip (-2.25,-2.0) rectangle (2.25,2.5);
\end{tikzpicture}}

\newcommand{\humanknot}{
\begin{tikzpicture}[thick,>=stealth,x={(3,0)},y={(0,2.5)},z={(0,0.05)},scale=0.875]
\def\p{{4,7,1,5,8,3,9,2,6,0,4}}
\def\q{{9,2,7,4,8,0,5,1,6,3,9}}
\foreach \s/\r in {0/solid,1/dashed,2/dashed,3/dashed,4/dashed,5/dashed,6/dashed,7/dashed,8/dashed,9/solid} { 
\draw[gray,\r,thick] (1,0,\s) \foreach \t in {5,10,...,180} {--({cos(\t)},{sin(\t)},\s)};}
\foreach \i in {5,7,1,9,3,6,8,2,4,0} {
\draw[>=triangle 90 cap,white,<->,line width = 4,shorten >=1, shorten <=1] 
({cos((\p[\i]+0.25)*36)},{sin((\p[\i]+0.25)*36)},\q[\i]) -- ({cos((\p[\i+1]+0.25)*36)},{sin((\p[\i+1]+0.25)*36)},\q[\i]); 
\draw[black,-,line width = 1,decoration={markings,mark=at position 0.4 with {\arrow{<}}},postaction={decorate}] 
({cos((\p[\i]+0.25)*36)},{sin((\p[\i]+0.25)*36)},\q[\i]) -- ({cos((\p[\i+1]+0.25)*36)},{sin((\p[\i+1]+0.25)*36)},\q[\i]); 
\filldraw ({cos((\p[\i]+0.25)*36)},{sin((\p[\i]+0.25)*36)},\q[\i]) circle[radius=0.01];
\draw[>=triangle 90 cap,white,<->,line width = 4,shorten >=1, shorten <=1] 
({cos((\p[\i+1]+0.25)*36)},{sin((\p[\i+1]+0.25)*36)},\q[\i]) -- ({cos((\p[\i+1]+0.25)*36)},{sin((\p[\i+1]+0.25)*36)},\q[\i+1]); 
\draw[black,-,line width = 1,decoration={markings,mark=at position 0.5 with {\arrow{<}}},postaction={decorate}] 
({cos((\p[\i+1]+0.25)*36)},{sin((\p[\i+1]+0.25)*36)},\q[\i]) -- ({cos((\p[\i+1]+0.25)*36)},{sin((\p[\i+1]+0.25)*36)},\q[\i+1]); 
\filldraw ({cos((\p[\i+1]+0.25)*36)},{sin((\p[\i+1]+0.25)*36)},\q[\i]) circle[radius=0.01];
\pgfmathsetmacro\pi{int(\p[\i])}
\draw ({1.2*cos((\p[\i]+0.25)*36)},{1.2*sin((\p[\i]+0.25)*36)},{9*int(\pi<5)}) node[fill=white] {${\mathlarger{\mathlarger{\mathlarger{\smiley}}}}_{\pgfmathprintnumber{\pi}}$};}
\foreach \s/\r in {0/solid,1/dashed,2/dashed,3/dashed,4/dashed,5/dashed,6/dashed,7/dashed,8/dashed,9/solid} { 
\draw[gray,\r,thick] (-1,0,\s) \foreach \t in {180,185,...,360} {--({cos(\t)},{sin(\t)},\s)};}
\filldraw ({cos((4.25)*36)},{sin((4.25)*36)},9) circle[radius=0.03];
\pgfresetboundingbox \clip (-1.25,-1.25) rectangle (1.25,1.75);
\end{tikzpicture}}

\newcommand{\knotgrid}{
\begin{tikzpicture}[thick,>=stealth]
\def\p{{4,7,1,5,8,3,9,2,6,0,4}}
\def\q{{9,2,7,4,8,0,5,1,6,3,9}}
\foreach \i in {5,7,1,9,3,6,8,2,4,0} {
\draw[>=triangle 90 cap,white,<->,line width = 4,shorten >=1, shorten <=1] 
({\p[\i]*0.5+0.5},{\q[\i]*0.5+0.5}) -- ({\p[\i+1]*0.5+0.5},{\q[\i]*0.5+0.5}); 
\draw[black,-,line width = 1,decoration={markings,mark=at position 0.45 with {\arrow{<}}},postaction=decorate] 
({\p[\i]*0.5+0.5},{\q[\i]*0.5+0.5}) -- ({\p[\i+1]*0.5+0.5},{\q[\i]*0.5+0.5});}
\foreach \i in {5,7,1,9,3,6,8,2,4,0} {
\draw[>=triangle 90 cap,white,<->,line width = 4,shorten >=1, shorten <=1] 
({\p[\i+1]*0.5+0.5},{\q[\i]*0.5+0.5}) -- ({\p[\i+1]*0.5+0.5},{\q[\i+1]*0.5+0.5});
\draw[black,-,line width = 1,decoration={markings,mark=at position 0.4 with {\arrow{<}}},postaction=decorate] 
({\p[\i+1]*0.5+0.5},{\q[\i]*0.5+0.5}) -- ({\p[\i+1]*0.5+0.5},{\q[\i+1]*0.5+0.5});
\filldraw ({\p[\i+1]*0.5+0.5},{\q[\i]*0.5+0.5}) circle[radius=0.02];
\filldraw ({\p[\i+1]*0.5+0.5},{\q[\i+1]*0.5+0.5}) circle[radius=0.02];}
\filldraw (2.5,5) circle[radius=0.05];
\draw[gray,->] (0,0) -- coordinate (x axis mid) (5.5,0);
\draw[gray,->] (0,0) -- coordinate (y axis mid) (0,5.5);
\foreach \x in {0,...,9} \draw (\x*0.5+0.5,1pt) -- (\x*0.5+0.5,-3pt) node[anchor=north] {\x}; 
\foreach \y in {0,...,9} \draw (1pt,\y*0.5+0.5) -- (-3pt,\y*0.5+0.5) node[anchor=east] {\y};
\pgfresetboundingbox \clip (-0.5,-0.5) rectangle (5.5,5.5);
\end{tikzpicture}}

\newcommand{\jumpknot}{
\begin{tikzpicture}[thick,>=stealth,x={(1,0)},y={(0,1)},z={(0.5,0)}, scale=2.4]
\draw[gray,dashed] (-1,0,0) arc (180:0:1cm and 0.5cm);
\draw[gray,dashed] (0,1,0) arc (90:-90:0.5cm and 1cm);
\draw[gray] (0,0) circle (1cm);
\shade[ball color=white,opacity=0.1] (0,0) circle (1cm);
\foreach \i/\x/\y/\z in {0/-0.192/-0.891/-0.187,1/-0.332/0.784/0.173,2/0.054/-0.530/-0.138,3/0.683/0.561/-0.007,4/0.522/-0.807/-0.653,5/-0.825/-0.514/0.267,6/-0.434/-0.067/-0.818,7/-0.809/-0.141/0.580,8/-0.715/0.805/0.483,9/0.614/0.311/0.257,10/0.218/-0.489/0.799,11/-0.914/0.555/0.556,12/0.577/0.358/-0.731,13/0.084/-0.292/0.921,14/0.597/0.878/-0.689,15/-0.650/0.075/-0.251,16/0.809/-0.011/0.060,17/-0.507/0.161/0.843,18/-0.505/0.027/-0.108,19/0.119/0.870/-0.314}
\node[fill, shading=ball, circle, ball color=gray, scale=0.5] (\i) at (\x,\y,\z) {};
\draw[white,line width=0.5pt,double=black,double distance=1.0pt] (10)--(11);
\draw[white,line width=0.5pt,double=black,double distance=1.0pt] (7)--(8);
\draw[white,line width=0.5pt,double=black,double distance=1.0pt] (9)--(10);
\draw[white,line width=0.5pt,double=black,double distance=1.0pt] (16)--(17);
\draw[white,line width=0.5pt,double=black,double distance=1.0pt] (8)--(9);
\draw[white,line width=0.5pt,double=black,double distance=1.0pt] (17)--(18);
\draw[white,line width=0.5pt,double=black,double distance=1.0pt] (13)--(14);
\draw[white,line width=0.5pt,double=black,double distance=1.0pt] (12)--(13);
\draw[white,line width=0.5pt,double=black,double distance=1.0pt] (1)--(2);
\draw[white,line width=0.5pt,double=black,double distance=1.0pt] (0)--(1);
\draw[white,line width=0.5pt,double=black,double distance=1.0pt] (2)--(3);
\draw[white,line width=0.5pt,double=black,double distance=1.0pt] (11)--(12);
\draw[white,line width=0.5pt,double=black,double distance=1.0pt] (15)--(16);
\draw[white,line width=0.5pt,double=black,double distance=1.0pt] (6)--(7);
\draw[white,line width=0.5pt,double=black,double distance=1.0pt] (4)--(5);
\draw[white,line width=0.5pt,double=black,double distance=1.0pt] (18)--(19);
\draw[white,line width=0.5pt,double=black,double distance=1.0pt] (19)--(0);
\draw[white,line width=0.5pt,double=black,double distance=1.0pt] (5)--(6);
\draw[white,line width=0.5pt,double=black,double distance=1.0pt] (3)--(4);
\draw[white,line width=0.5pt,double=black,double distance=1.0pt] (14)--(15);
\draw[gray] (-1,0,0) arc (180:360:1cm and 0.5cm);
\draw[gray] (0,1,0) arc (90:270:0.5cm and 1cm);
\pgfresetboundingbox \clip (-1.0,-1.0) rectangle (1.0,1.0);
\end{tikzpicture}}

\newcommand{\jumplink}{
\begin{tikzpicture}[thick,>=stealth,x={(0.6,0)},y={(0,0.6)},z={(0.3,0.2)}, scale=3.3]
\draw[gray,dashed] (-1,-1,-1) -- (-1,-1,+1);
\draw[gray] (-1,+1,-1) -- (-1,+1,+1);
\draw[gray] (+1,-1,-1) -- (+1,-1,+1);
\draw[gray] (-1,-1,-1) -- (-1,+1,-1);
\draw[gray,dashed] (-1,-1,+1) -- (-1,+1,+1);
\draw[gray] (+1,-1,+1) -- (+1,+1,+1);
\draw[gray] (-1,-1,-1) -- (+1,-1,-1);
\draw[gray,dashed] (-1,-1,+1) -- (+1,-1,+1);
\draw[gray] (-1,+1,+1) -- (+1,+1,+1);
\foreach \i/\x/\y/\z in {00/-0.635/0.886/0.722,01/0.708/-0.846/0.593,02/0.709/0.578/0.907,03/0.415/-0.829/0.203,04/-0.935/-0.581/0.282,05/0.513/-0.115/-0.717,06/0.675/0.925/0.616,07/-0.983/-0.233/-0.553,10/0.941/-0.790/0.658,11/-0.508/-0.087/-0.294,12/0.371/0.399/0.091,13/-0.325/-0.820/-0.796,14/0.121/-0.594/0.157,15/-0.799/0.538/0.137,16/-0.908/0.497/-0.370,17/0.901/-0.172/-0.251}
\node[fill, shading=ball, circle, ball color=gray, scale=0.5] (\i) at (\x,\y,\z) {};
\draw[white,line width=0.5pt,double=lightgray,double distance=1.5pt] (00)--(01);
\draw[white,line width=0.5pt,double=lightgray,double distance=1.5pt] (03)--(04);
\draw[white,line width=0.5pt,double=black,double distance=1.5pt] (17)--(10);
\draw[white,line width=0.5pt,double=black,double distance=1.5pt] (10)--(11);
\draw[white,line width=0.5pt,double=black,double distance=1.5pt] (14)--(15);
\draw[white,line width=0.5pt,double=lightgray,double distance=1.5pt] (07)--(00);
\draw[white,line width=0.5pt,double=lightgray,double distance=1.5pt] (06)--(07);
\draw[white,line width=0.5pt,double=lightgray,double distance=1.5pt] (05)--(06);
\draw[white,line width=0.5pt,double=black,double distance=1.5pt] (11)--(12);
\draw[white,line width=0.5pt,double=black,double distance=1.5pt] (15)--(16);
\draw[white,line width=0.5pt,double=lightgray,double distance=1.5pt] (01)--(02);
\draw[white,line width=0.5pt,double=lightgray,double distance=1.5pt] (04)--(05);
\draw[white,line width=0.5pt,double=black,double distance=1.5pt] (16)--(17);
\draw[white,line width=0.5pt,double=black,double distance=1.5pt] (13)--(14);
\draw[white,line width=0.5pt,double=lightgray,double distance=1.5pt] (02)--(03);
\draw[white,line width=0.5pt,double=black,double distance=1.5pt] (12)--(13);
\draw[gray] (+1,+1,-1) -- (+1,+1,+1);
\draw[gray] (+1,-1,-1) -- (+1,+1,-1);
\draw[gray] (-1,+1,-1) -- (+1,+1,-1);
\pgfresetboundingbox \clip (-1.5,-1.25) rectangle (1.5,1.25);
\end{tikzpicture}}

\newcommand{\curvegraph}{
\tikzstyle{every node}=[circle, draw, fill=black!50, minimum width=6pt, inner sep=0pt]
\begin{tikzpicture}[thick,scale=0.75,rotate=270]
\draw [black,line width=2.0pt] (2,5) to[out=135,in=90] (0,4);
\draw [black,line width=2.0pt] (4,6) to[out=270,in=90] (1,5);
\draw [black,line width=2.0pt] (2,3) to[out=0,in=270] (4,5);
\draw [black,line width=2.0pt] (4,6) to[out=90,in=90] (3,6);
\draw [black,line width=2.0pt] (3,6) to[out=270,in=90] (4,5);
\draw [black,line width=2.0pt] (2,4) to[out=180,in=270] (1,5);
\draw [black,line width=2.0pt] (5,3) to[out=90,in=0] (4,5);
\draw [black,line width=2.0pt] (4,3) to[out=90,in=0] (2,4);
\draw [black,line width=2.0pt] (3,4) to[out=270,in=90] (3,2);
\draw [black,line width=2.0pt] (3,4) to[out=90,in=180] (4,5);
\draw [black,line width=2.0pt] (2,0) to[out=0,in=270] (3,2);
\draw [black,line width=2.0pt] (3,1) to[out=0,in=270] (5,3);
\draw [black,line width=2.0pt] (2,2) to[out=180,in=90] (1,1);
\draw [black,line width=2.0pt] (3,1) to[out=180,in=270] (1,3);
\draw [black,line width=2.0pt] (2,2) to[out=0,in=180] (3,2);
\draw [black,line width=2.0pt] (2,5) to[out=315,in=90] (1,3);
\draw [black,line width=2.0pt] (4,3) to[out=270,in=0] (3,2);
\draw [black,line width=2.0pt] (2,0) to[out=180,in=270] (1,1);
\draw [black,line width=2.0pt] (2,3) to[out=180,in=270] (0,4);
\node at(2.9,1.0){};
\node at(1.4,1.8){};
\node at(3.0,2.0){};
\node at(1.0,3.1){};
\node at(3.0,3.3){};
\node at(3.0,3.9){};
\node at(3.6,3.75){};
\node at(1.7,4.1){};
\node at(4.0,5.0){};
\node at(1.1,5.2){};
\node at(3.5,5.5){};
\pgfresetboundingbox \clip (-0.5,0) rectangle (5.5,7);
\end{tikzpicture}}

\newcommand{\curveknot}{
\begin{tikzpicture}[thick,scale=0.75,rotate=270]
\draw [white,line width=3.0pt,double=black,double distance=2.0pt] (2,5) to[out=135,in=90] (0,4);
\draw [white,line width=3.0pt,double=black,double distance=2.0pt] (4,6) to[out=270,in=90] (1,5);
\draw [white,line width=3.0pt,double=black,double distance=2.0pt] (2,3) to[out=0,in=270] (4,5);
\draw [white,line width=3.0pt,double=black,double distance=2.0pt] (4,6) to[out=90,in=90] (3,6);
\draw [white,line width=3.0pt,double=black,double distance=2.0pt] (3,6) to[out=270,in=90] (4,5);
\draw [white,line width=3.0pt,double=black,double distance=2.0pt] (2,4) to[out=180,in=270] (1,5);
\draw [white,line width=3.0pt,double=black,double distance=2.0pt] (5,3) to[out=90,in=0] (4,5);
\draw [white,line width=3.0pt,double=black,double distance=2.0pt] (4,3) to[out=90,in=0] (2,4);
\draw [white,line width=3.0pt,double=black,double distance=2.0pt] (3,4) to[out=270,in=90] (3,2);
\draw [white,line width=3.0pt,double=black,double distance=2.0pt] (3,4) to[out=90,in=180] (4,5);
\draw [white,line width=3.0pt,double=black,double distance=2.0pt] (2,0) to[out=0,in=270] (3,2);
\draw [white,line width=3.0pt,double=black,double distance=2.0pt] (3,1) to[out=0,in=270] (5,3);
\draw [white,line width=3.0pt,double=black,double distance=2.0pt] (2,2) to[out=180,in=90] (1,1);
\draw [white,line width=3.0pt,double=black,double distance=2.0pt] (3,1) to[out=180,in=270] (1,3);
\draw [white,line width=3.0pt,double=black,double distance=2.0pt] (2,2) to[out=0,in=180] (3,2);
\draw [white,line width=3.0pt,double=black,double distance=2.0pt] (2,5) to[out=315,in=90] (1,3);
\draw [white,line width=3.0pt,double=black,double distance=2.0pt] (4,3) to[out=270,in=0] (3,2);
\draw [white,line width=3.0pt,double=black,double distance=2.0pt] (2,0) to[out=180,in=270] (1,1);
\draw [white,line width=3.0pt,double=black,double distance=2.0pt] (2,3) to[out=180,in=270] (0,4);
\draw[black,-,line width = 2,decoration={markings,mark=at position 0.45 with {\arrow{<}}},postaction=decorate] (5,3) to[out=90,in=0] (4,5);
\pgfresetboundingbox \clip (-0.5,0) rectangle (5.5,7);
\end{tikzpicture}}

\newcommand{\squarecross}[2]
{\tikzstyle{every node}=[circle, draw, fill=black!50, inner sep=0pt, minimum width=3pt]
\begin{tikzpicture}[thick,baseline=-4pt,scale=0.75]
\node at(-0.2,0.2)(0){}; \node at(-0.2,-0.2)(1){};
\node at (0.2,0.2)(2){}; \node at (0.2,-0.2)(3){};
\draw[white,line width=1.5pt,double=black,double distance=1.0pt] { #1 };
\draw[white,line width=1.5pt,double=black,double distance=1.0pt] { #2 };
\end{tikzpicture}}

\newcommand{\braidknot}{
\begin{tikzpicture}[scale=0.5]
\foreach \x/\ya/\yb/\yc/\za/\zb/\zc in {
0/0/1/2/1/0/2,
1/0/1/2/1/0/2,
2/0/2/1/0/1/2,
3/0/2/1/0/1/2,
4/0/1/2/1/0/2,
5/0/2/1/0/1/2,
6/1/0/2/0/1/2,
7/0/1/2/0/2/1,
8/0/1/2/1/0/2,
9/0/1/2/1/0/2}{
\draw[white,line width=3.0pt,double=black,double distance=2.0pt] (\x,\yc) to[out=0,in=180] ({\x+1},\zc);
\draw[white,line width=3.0pt,double=black,double distance=2.0pt] (\x,\yb) to[out=0,in=180] ({\x+1},\zb);
\draw[white,line width=3.0pt,double=black,double distance=2.0pt] (\x,\ya) to[out=0,in=180] ({\x+1},\za);}
\foreach \i in {0,1,2}{
\draw[lightgray,line width=2.0pt] (10,\i) to[out=0,in=0] ({10-0.1*\i},{3.2-0.2*\i}) -- ({0+0.1*\i},{3.2-0.2*\i}) to[out=180,in=180] (0,\i);}
\pgfresetboundingbox \clip (-1.5,0) rectangle (11.5,4.5);
\end{tikzpicture}}

\newcommand{\bridgeknot}{
\begin{tikzpicture}[scale=0.5]
\foreach \x/\ya/\yb/\yc/\yd/\za/\zb/\zc/\zd in {
0/1/0/2/3/0/1/2/3,
1/0/1/2/3/0/1/3/2,
2/0/2/1/3/0/1/2/3,
3/0/2/1/3/0/1/2/3,
4/0/1/2/3/0/1/3/2,
5/0/1/3/2/0/1/2/3,
6/0/1/2/3/1/0/2/3,
7/0/1/2/3/0/2/1/3,
8/1/0/2/3/0/1/2/3}{
\draw[white,line width=3.0pt,double=black,double distance=2.0pt] (\x,\zd) to[out=0,in=180] ({\x+1},\yd);
\draw[white,line width=3.0pt,double=black,double distance=2.0pt] (\x,\zc) to[out=0,in=180] ({\x+1},\yc);
\draw[white,line width=3.0pt,double=black,double distance=2.0pt] (\x,\zb) to[out=0,in=180] ({\x+1},\yb);
\draw[white,line width=3.0pt,double=black,double distance=2.0pt] (\x,\za) to[out=0,in=180] ({\x+1},\ya);}
\draw[lightgray,line width=2.0pt] (9,0) -- (9.4,0) to[out=0,in=0] (9.4,1) --(9,1);
\draw[lightgray,line width=2.0pt] (9,3) -- (9.4,3) to[out=0,in=0] (9.4,2) --(9,2);
\draw[lightgray,line width=2.0pt] (0,0) -- (-0.4,0) to[out=180,in=180] (-0.4,1) -- (0,1);
\draw[lightgray,line width=2.0pt] (0,3) -- (-0.4,3) to[out=180,in=180] (-0.4,2) -- (0,2);
\pgfresetboundingbox \clip (-1,0) rectangle (10,4.5);
\end{tikzpicture}}

\newcommand{\lissajouscurve}{
\begin{tikzpicture}[scale=1]
\foreach \t in {0,0.001,...,1}
{\draw[black,line width=1.0pt] ({cos(4*180*(\t)+30)},{cos(6*180*(\t))}) -- ({cos(4*180*(\t+0.001)+30)},{cos(6*180*(\t+0.001))});}
\pgfresetboundingbox \clip (-1,-1) rectangle (1,1.4);
\end{tikzpicture}}

\newcommand{\lissajousknot}{
\begin{tikzpicture}[scale=1]
\foreach \i in {1, 12, 0, 14, 23, 21, 16, 19, 4, 10, 13, 22, 5, 8, 17, 9, 11, 6, 20, 2, 18, 3, 7, 15} {
\draw[white,line width=3pt,double=black,double distance=3pt] 
({abs(Mod(2*\i+1,6*4)/4-3)},{abs(Mod(2*\i-1,4*4)/4-2)}) -- 
({abs(Mod(2*\i+2,6*4)/4-3)},{abs(Mod(2*\i+0,4*4)/4-2)}) -- 
({abs(Mod(2*\i+3,6*4)/4-3)},{abs(Mod(2*\i+1,4*4)/4-2)}); }
\pgfresetboundingbox \clip (0,0) rectangle (3,2.4);
\end{tikzpicture}}

\newcommand{\harmoniccurve}{
\begin{tikzpicture}[scale=1]
\foreach \t in {0,0.001,...,1}
{\draw[black,line width=1.0pt] ({cos(5*180*(\t))},{cos(9*180*(\t))}) -- ({cos(5*180*(\t+0.001))},{cos(9*180*(\t+0.001))});}
\pgfresetboundingbox \clip (-1,-1) rectangle (1,1.4);
\end{tikzpicture}}

\newcommand{\harmonicknot}{
\begin{tikzpicture}[scale=0.8]
\foreach \i in {39, 5, 24, 9, 22, 8, 40, 31, 16, 3, 34, 33, 44, 4, 35, 30, 45, 23, 21, 13, 19, 0, 38, 42, 14, 27, 10, 29, 28, 41, 12, 18, 17, 32, 20, 25, 26, 43, 1, 6, 15, 11, 2, 7, 36, 37} {
\draw[white,line width=3pt,double=black,double distance=3pt] 
({abs(Mod(2*\i-1,9*4)/4-4.5)},{abs(Mod(2*\i-1,5*4)/4-2.5)}) -- 
({abs(Mod(2*\i+0,9*4)/4-4.5)},{abs(Mod(2*\i+0,5*4)/4-2.5)}) -- 
({abs(Mod(2*\i+1,9*4)/4-4.5)},{abs(Mod(2*\i+1,5*4)/4-2.5)}); }
\pgfresetboundingbox \clip (0,0) rectangle (4.5,3);
\end{tikzpicture}}

\newcommand{\starknot}{
\begin{tikzpicture}[thick,scale=0.5625]
\foreach \x in {0,1,...,8} {\draw[white,-,line width = 3] (\x*280-70:4) -- (\x*280+90:4); \draw[black] (\x*280-70:4) -- (\x*280+90:4); }
\foreach \x in {0,1,...,8} { \draw[->] (\x*160-70:4) -- node[auto,swap,pos=0.5]{\small} +(\x*160+100:1); \draw (\x*280-70:4) -- +(\x*280+120:1); }
\pgfresetboundingbox \clip (-4,-4) rectangle (4,4);
\end{tikzpicture}}

\newcommand{\starbraid}{
\tikz[thick,rotate=180]{
\foreach \x/\ya/\yb/\yc/\yd/\za/\zb/\zc/\zd/\ua/\ub/\uc/\ud/\va/\vb/\vc/\vd in {
 0/1/0/2/3/0/1/3/2/0/2/1/3/0/1/2/3,
 2/0/1/3/2/1/0/2/3/0/2/1/3/0/1/2/3,
 4/1/0/2/3/0/1/3/2/0/2/1/3/0/1/2/3,
 6/0/1/3/2/1/0/2/3/0/2/1/3/0/1/2/3,
 8/1/0/2/3/0/1/3/2/0/2/1/3/0/1/2/3,
10/0/1/2/3/1/0/3/2/0/1/2/3/0/2/1/3,
12/0/1/2/3/1/0/3/2/0/2/1/3/0/1/2/3,
14/1/0/2/3/0/1/3/2/0/2/1/3/0/1/2/3,
16/0/1/3/2/1/0/2/3/0/2/1/3/0/1/2/3}{
\draw[white,line width=2.0pt,double=black,double distance=1.0pt] (\x*20:0.3*\zd+1.1) -- (\x*20+20:0.3*\yd+1.1);
\draw[white,line width=2.0pt,double=black,double distance=1.0pt] (\x*20:0.3*\zc+1.1) -- (\x*20+20:0.3*\yc+1.1);
\draw[white,line width=2.0pt,double=black,double distance=1.0pt] (\x*20:0.3*\zb+1.1) -- (\x*20+20:0.3*\yb+1.1);
\draw[white,line width=2.0pt,double=black,double distance=1.0pt] (\x*20:0.3*\za+1.1) -- (\x*20+20:0.3*\ya+1.1);
\draw[black,line width=1.0pt,decoration={markings,
   mark=at position 0.5 with {\arrow{>}}},postaction={decorate}] (\x*20+20:0.3*\vd+1.1) -- (\x*20+40:0.3*\ud+1.1);
\draw[white,line width=2.0pt,double=black,double distance=1.0pt] (\x*20+20:0.3*\vc+1.1) -- (\x*20+40:0.3*\uc+1.1);
\draw[white,line width=2.0pt,double=black,double distance=1.0pt] (\x*20+20:0.3*\vb+1.1) -- (\x*20+40:0.3*\ub+1.1);
\draw[white,line width=2.0pt,double=black,double distance=1.0pt] (\x*20+20:0.3*\va+1.1) -- (\x*20+40:0.3*\ua+1.1);}
\pgfresetboundingbox \clip (-2.25,-2.25) rectangle (2.25,2.25);}}

\title{
Models of Random Knots}

\author{
Chaim Even-Zohar
\footnote{
Department of Mathematics, University of California, Davis, California 95616. chaim@math.ucdavis.edu
}}

\date{}

\begin{document}

\maketitle

\begin{abstract}
The study of knots and links from a probabilistic viewpoint provides insight into the behavior of ``typical'' knots, and opens avenues for new constructions of knots and other topological objects with interesting properties. The knotting of random curves arises also in applications to the natural sciences, such as in the context of the structure of polymers. 

We present here several known and new randomized models of knots and links. We review the main known results on the knot distribution in each model. We discuss the nature of these models and the properties of the knots they produce.

Of particular interest to us are finite type invariants of random knots, and the recently studied Petaluma model. We report on rigorous results and numerical experiments concerning the asymptotic distribution of such knot invariants. Our approach raises questions of universality and classification of the various random knot models.   

\smallskip \noindent \textbf{MSC} 57M25, 60B05
\end{abstract}

\section{Knots}

There is an increasing interest in random knots by both topologists and probabilists, as well as researchers from other disciplines. Our aim in this survey article is to provide an accessible overview of the many different approaches to this topic. 

We start with a very brief introduction to knot theory, and in Section~\ref{randomness-section} we describe the motivation to introduce randomness into this field. The various models are surveyed in Section~\ref{models-section}, and some specific aspects are further discussed in Section~\ref{new-section}. Some thoughts and open problems are presented in Section~\ref{discussion}.

\medskip

Intuitively, a knot is a simple closed curve in the three dimensional space, considered up to continuous deformations without self-crossing. More formally, a \emph{knot} is a smoothly embedded oriented circle $S^1 \hookrightarrow \R^3$, with the equivalence relation of ambient isotopies of $\R^3$. A \emph{link} is a disjoint union of several such embedded circles, called \emph{components}, with the same equivalence. An alternative definition uses polygonal paths without the smoothness condition. There are several good general introductions to knot theory such as~\cite{adams1994knot} or~\cite{lickorish1997introduction}.

\begin{figure}[t]
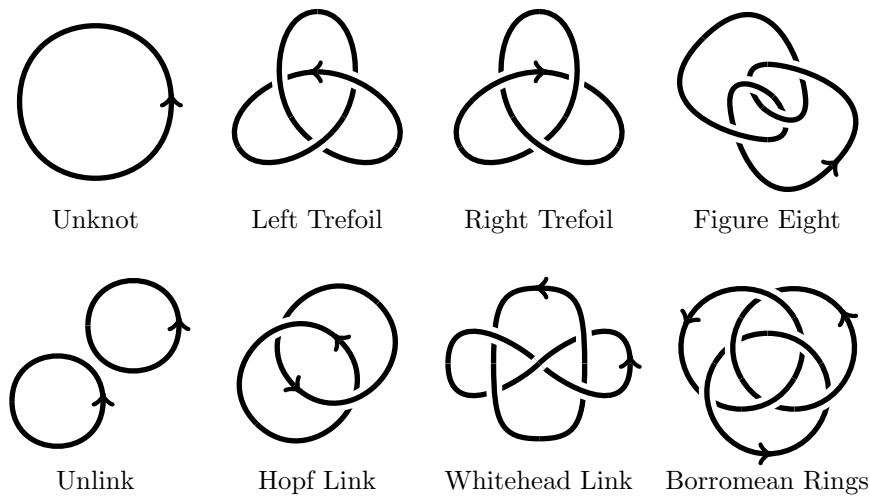

\begin{center}
\knottable
\caption{Selected knot and link diagrams.}
\label{diagrams}
\end{center}
\end{figure}

Knots and links can also be described via planar \emph{diagrams}, which are their generic projections to~$\R^2$. The projection is injective except for a finite number of traverse double points. Each such \emph{crossing} point is decorated to indicate which preimage is over and which is under, with respect to the direction of the projection. See Figure~\ref{diagrams} for diagrams of some well-known knots and links. 

The set of nonequivalent knots is infinite, without much structure and organization. Some order arise from the operation of \emph{connected sum} of knots, $\trefoil \# \trefoil = \granny$ for example.  A~theorem by Schubert~\cite[Chapter 2]{lickorish1997introduction} states that every knot can be uniquely decomposed as a connected sum of \emph{prime} knots, which are knot that cannot be decomposed further. However, there are infinitely many nonequivalent prime knots as well. 

A problem that motivated much of the developments in knot theory since its early days was finding and tabulating all prime knots that can be represented by diagrams with a small number of crossings. As of today, knot tables with up to 16 crossings have been compiled~\cite{hoste1998first}. This classification mission called for tools for telling whether or not two given knots are equivalent, even though represented differently.

By the classical Reidemeister theorem~\cite[Chapter 1]{lickorish1997introduction}, two diagrams define equivalent knots if and only if one can be transformed into the other by a sequence of local \emph{moves} of three types: (I)~twisting the curve at some point~\mbox{$\;\RMIa \leftrightarrow \RMIb\;$}, (II)~sliding one part of the curve under an adjacent part~\mbox{$\;\RMIIa \leftrightarrow \RMIIb\;$}, or (III)~sliding under an adjacent crossing~\mbox{$\;\RMIIIa \leftrightarrow \RMIIIb\;$}. 

As for the complementary purpose of distinguishing one knot from another, a wide variety of \emph{knot invariants} were defined. Here one constructs a well-defined function from the set of all knots to any other set, that attains different values for the two knots in question. Either representation, by diagrams or by curves in $\R^3$, may be used to define invariants, as long as one shows that it respects equivalence. In a broader perspective, knot invariants may be viewed as tools to classify knots and understand their properties.

\medskip

We mention some important knot invariants. The \emph{crossing number} $c(K)$ is the least number of crossing points in a diagram of a knot $K$. The \emph{genus} $g(K)$ is the least genus of an embedded oriented compact surface with boundary $K$. Several other invariants, such as the \emph{bridge number}, \emph{unknotting number}, and \emph{stick number}~\cite{adams1994knot}, are similarly defined by taking the minimum value of some complexity measure over certain descriptions of the knot.  

It is conjectured that knots can be fully classified by Vassiliev's \emph{finite type invariants}~\cite{vassiliev1988cohomology, bar1995vassiliev, chmutov2012introduction}. See Section~\ref{fti-def} for a definition. This infinite collection of numerical invariants includes Gauss's \emph{linking number} and the \emph{Casson invariant}, coefficients of the \emph{Alexander--Conway polynomial}, the modified \emph{Jones polynomial}, and the \emph{Kontsevich integral}.  

Other invariants are defined via properties of the knot complement, such as its fundamental group $\pi_1(\R^3 \setminus K)$, the \emph{knot group} of $K$. A knot is called \emph{hyperbolic} if its complement can be given a metric of constant negative curvature. In this case $\text{Vol}(S^3 \setminus K)$, the \emph{hyperbolic volume} of $K$, is a well-defined and useful knot invariant~\cite{thurston1978geometry}.

\medskip

\section{Randomness}\label{randomness-section}

There are several motivations to study randomized knot models. They emerge from different perspectives. Below we mention several aspects and applications of knot theory where it is natural to adopt a probabilistic point of view.

\paragraph{Study Typical Knots} 

As metioned, the space of knots is infinite and poorly structured. Usually, the particular examples of knots one considers are either very simple with only a few crossings, or they are explicit constructions of knots of quite specific forms. These can be members of well-known families such as torus knots, pretzel knots, and rational knots, or ad hoc constructions tailored for the problem under investigation. 

Similarly, often one considers knots of some particular type, such as alternating or hyperbolic. Do these classes represent the general case, and if so in what sense?

It is natural and desirable to understand what typical knots are like and what properties they tend to have. 

\medskip

We specify a probability distribution over knots in search of a framework to investigate such questions. Often we consider a sequence of such distributions, supported on increasingly larger sets of knots. These distributions may be defined via random planar diagrams or random curves in $\R^3$, but ultimately only the resulting knot type is considered. 

Rather than focus on particular constructions and classes, we ask what knot properties hold with high probability. Knot invariants become random variables on the probability space, and we study their distribution and interrelations. 

It is not apriori clear which distributions, or \emph{models of random knots}, are worth studying. It is reasonable to require that every knot have positive probability. We also do not want the measure to be highly concentrated on some overly specific class of knots. At present it remains debatable how good any concrete distribution that one suggests is.

\paragraph{Probabilistic Existence Proofs} 

A more definite goal of studying random knot models is the application of the probabilistic method in knot theory. The basic idea is to prove the existence of certain objects by showing that in some random model they occur with positive probability. This influential methodology has yielded many unexpected results in combinatorics and other fields~\cite{alon2000probabilistic}. In many cases, the existence of objects with some given properties can be established using probabilistic methods, while finding matching explicit constructions remains elusive. 

To illustrate this idea, consider the Jones polynomial $V_K(t) \in \Z\left[t,t^{-1}\right]$. The discovery of this important knot invariant in 1984 was hailed as a breakthrough in the field~\cite[Chapter 3]{lickorish1997introduction}. By definition $V_{\text{unknot}}(t)=1$, and it is unknown whether there exists a non-trivial knot $K$ for which $V_K(t)=1$~\cite{jones2000ten}. It is believed though, that if such knots exist, then they are plentiful. If so, it is reasonable to expect that in some random model it should be possible to prove the probability of this trivial Jones polynomial is strictly larger than that of the unknot.

For this approach to work we clearly need a random model that allows us to estimate the probability of the relevant events and the distributions of the invariants at hand. 

Random knot models come handy also in computer experiments, where one non-exhaustively searches for a specific exemplar to demonstrate some properties, thus providing explicit examples and counterexamples in a more direct way.

\paragraph{Knots in Nature}

The occurrence of knots and links in the natural sciences has been a fruitful source for several studies of randomized knot models. 

\medskip

Most prominently, biologists are interested in the three-dimensional shape of proteins, DNA and RNA molecules. Their geometric and topological features affect their functionality in a variety of biological processes, such as protein folding and DNA replication and transcription. Physicists and chemists look into the formation of entanglements in polymeric materials. The topological structure of such substances is reflected macroscopically in its features, such as elasticity, viscosity, diffusion rate, and purity of crystallization. There is plenty of literature on the modeling of knots in thread-like molecules. Find some expositions and surveys in~\cite{wasserman1986biochemical, vologodskii1992topology, sumners1992knot, sumners1995lifting, grosberg1997giant, bates2005dna, orlandini2007statistical, mcleish2008tangled, fenlon2008open, buck2009dna, de2009lectures, micheletti2011polymers, lim2015molecular}.

Numerous numerical simulations and experiments have been preformed to investigate the topological properties of such filamentary molecules. These involved the invention of several mathematical models that produce random paths in $\R^3$ to simulate the conformation of molecules in natural environments. In particular, such a model defines a distribution over knot types, often parametrized by the length of the path.

Naturally, these models are designated to emulate natural features and processes, with different degrees of simplification. Most often they incorporate physical constraints such as non-zero thickness, self interaction, restricted bending, and spatial confinement. Additionally, this line of research calls for random models that can be easily sampled in numerical studies.

\begin{figure}[bt] \center
\includegraphics[width=0.45\columnwidth]{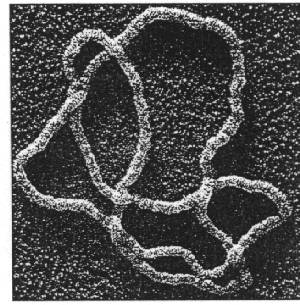}
\;\;\;
\includegraphics[width=0.45\columnwidth]{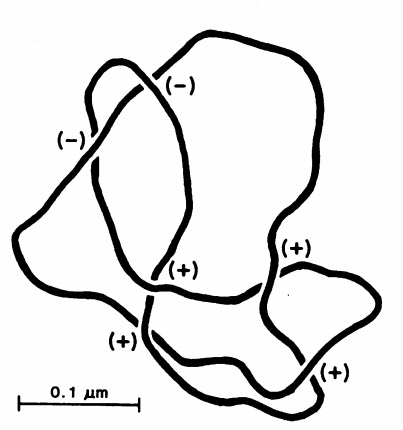}
\caption{Knotted DNA. Figure is courtesy of Wasserman et al.~\cite{wasserman1985discovery}}
\label{dna}
\end{figure}	

\medskip

The study of knotted structures in three-dimensional fields dates back to early fluid dynamics and Kelvin's vortex atom hypothesis~\cite{helmholtz1867lxiii, kelvin1867vortex}. Knots and links are formed in a three-dimensional flow $\vec{u}:\R^3 \to \R^3$ by the vortex lines that follow $\nabla \times \vec{u}$, or similarly by the nodal set $\psi=0$ of a wavefunction \mbox{$\psi:\R^3 \to \C$}. 

In current research, such knotting phenomena are theoretically analyzed, numerically simulated, and experimentally created or identified in various physical systems. To mention some examples: knotted vortices in classical fluid flow~\cite{kleckner2013creation} and in superfluids~\cite{hall2016tying, kleckner2016superfluid}, optical vortices in laser beams~\cite{dennis2010isolated}, magnetic fields in plasma~\cite{berger1999introduction}, superposition of states in quantum mechanics~\cite{berry2001knotted}, and also nonlinear waves in biological and chemical excitable media~\cite{winfree1984organizing}.

It seems that the generation of such knotted fields is often dominated by random factors, and it would be interesting to investigate what knots and links are likely to occur in such circumstances. Indeed, a recent work~\cite{taylor2016vortex} simulates random quantum wavefunctions in different potentials, and study the complexity of the vortex knots that show up.

\medskip

Finally, knots form at random in many objects of everyday practice, from extension cords, ropes, and garden hoses~\cite{raymer2007spontaneous} to umbilical cords~\cite{goriely2005knotted, hershkovitz2001risk} and eels~\cite{zintzen2011hagfish}.

\paragraph{Computational Aspects}

The study of random knot models is also motivated by the important role of randomness in the design and analysis of algorithms and in computational complexity theory~\cite{mitzenmacher2005probability}. 

It is a central computational challenge in knot theory to determine how hard it is to detect unknots, and more generally to decide the equivalence of two given knots~\cite{haken1961theorie,hass1999computational,kuperberg2014knottedness,lackenby2016efficient}. Specifically it is interesting to bound the number of Reidemeister moves that yield the equivalence of two representations~\cite{hass2010unknot,lackenby2015polynomial,coward2014upper}. It is generally believed that some of these problems are hard, and consequently cryptosystems were proposed that are based on such problems~\cite{farhi2012quantum}. To this end it is necessary to know the complexity of typical instances of problems. Random knot models are clearly needed in such pursuits.

The computation of various invariants also leads to interesting complexity problems. Hardness results are known for the knot genus~\cite{agol2006computational,lackenby2016efficient} and for the Jones polynomial~\cite{jaeger1990computational, aharonov2009polynomial, kuperberg2009hard}. Other invariants such as the Alexander--Conway polynomial and finite type invariants are computable in polynomial time~\cite{alexander1928topological,bar1995polynomial,chmutov2012introduction}. Many such algorithms are implemented in software packages, such as~\emph{SnapPy}~\cite{snappy}, \emph{KnotTheory}~\cite{knottheory} and \emph{KnotScape}~\cite{knotscape}. These are used in practice for the compilation of knot databases and are important tools in research and applications~\cite{knotatlas, knotinfo, hoste1998first}. Random knot models could serve as the basis for average-case analysis of such algorithms.

\paragraph{Random $3$-Manifolds}

The probabilistic method has had a great success in many areas. The study of random knots can be viewed as part of a broader research effort to apply this approach to the study of geometric and topological objects. 
In recent years, there have been interesting developments in the study of random simplicial complexes~\cite{linial2006homological, adler2010persistent, kahle2016random}, random groups~\cite{gromov2003random, ollivier2005january}, random manifolds~\cite{brooks2004random, dunfield2006finite, pippenger2006topological, farber2008betti}, and more. 

In particular, several models for random $3$-manifolds have been presented and studied in the past decade~\cite{dunfield2006finite, lutz2008combinatorial, maher2010random, kowalskicomplexity, dunfield2011quantum, maher2011random, maher2012exponential, lubotzky2016random, rivin2014statistics}. Since every closed orientable $3$-manifold can be generated by performing Dehn surgeries on links in $S^3$ \cite[Chapter 12]{lickorish1997introduction}, models for random links give rise to random $3$-manifolds whose properties are interesting to study~\cite{workinprogress}.

In another direction, random knot models may extend to knotted 2-spheres or other surfaces in a 4-sphere, and further to randomly embedded manifolds in higher dimensions~\cite{soteros2012knotted, atapour2015counting}.

~

\section{Models}\label{models-section}

Before listing some specific models, a few words on the general framework. A \emph{random knot model} is a distribution over the set of all knots, which we represent by a random variable $K$. We usually consider a sequence $K_n$ of such distributions, where $n \in \N$ naturally appears in the construction. This parameter $n$ can often be viewed as a complexity measure of the typical resulting knots. All unspecified asymptotic statements that we make here are w.r.t.~$n\to\infty$.

Variations abound: We also encounter some multi-parameter constructions and some models that yield random links of any number of components, or focus on some subclass such as prime or alternating knots.

\paragraph{Self-Avoiding Grid Walk} \label{walk-model}

\begin{figure}[bht]
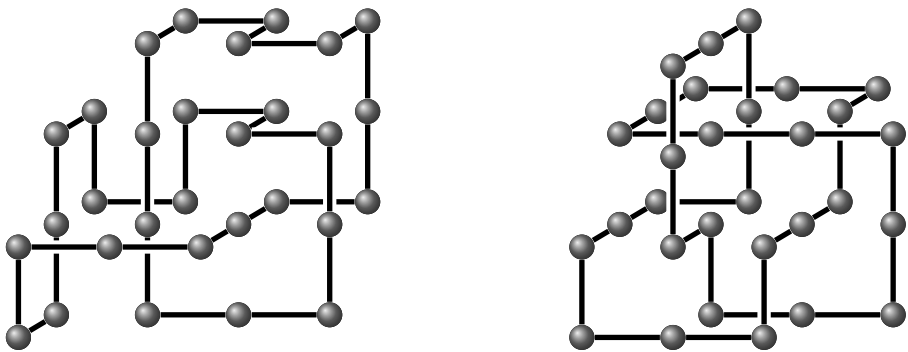

\begin{center}
\knotwalk
\caption{The trefoil and figure-eight knots as 30-step walks on $\Z^3$.}
\label{walks}
\end{center}
\end{figure}

As usual, a walk on the three-dimensional lattice $\Z^3$ is a sequence $\{X_0,\dots,X_n\}$ such that $X_0 = (0,0,0)$ and $(X_{i+1} - X_i) \in \{(\pm1,0,0),(0,\pm1,0),(0,0,\pm1)\}$. Consider $n$-step walks that are \emph{closed}, with $X_n=X_0$, and \emph{self-avoiding}, so that $X_i \neq X_j$ for any other pair of points. Connecting the points of such a walk yields an $n$-segment polygonal path, that represents a knot. See Figure~\ref{walks} for two examples. 

Random self-avoiding walks (SAW) on $\Z^3$ were suggested as a model for polymeric molecules, and their knotting properties have been studied over the past decades in varying degrees of rigor. In the \emph{grid walk model} a random knot $K_n$ is obtained by sampling from the uniform distribution of all closed self-avoiding $n$-step walks. Every knot appears in this model for $n$ large enough.

It was conjectured by Delbruck~\cite{delbruck1961knotting} that $K_n$ is knotted with high probability. This was observed in numerical simulations~\cite{crippen1974topology,frank1975statistical}, and proved by Sumners and Whittington~\cite{sumners1988knots} and by Pippenger~\cite{pippenger1989knots}. Using Kesten's pattern theorem~\cite{kesten1963number}, they showed that the unknot appears with exponentially small probability in $n$. Moreover, every prime knot appears in the decomposition of $K_n$ with multiplicity $\Theta(n)$, except for an exponentially small probability~\cite{soteros1992entanglement}. 

Let $K_n'$ be a uniformly random connected component of $K_n$, conditioned on $K_n$ being knotted. Note that $K_n'$ is a natural model for \emph{random prime knots}, and it is suggestive that $K_n'$ converges in distribution, and yields a random model for all prime knots.

It is of interest to study self-avoiding walks and the resulting knots on other lattices~\cite{van1990knot}. Also extensions to random 2-component links, and the effect of confinement within a box or a tube, are considered and analyzed~\cite{orlandini1994random, soteros1999linking, atapour2010linking}. Madras and Slade's book~\cite{madras2013self} offers a rigorous analysis of Monte Carlo sampling methods for self avoiding walks.

\paragraph{Polygonal Walks} \label{polygonal-model}

In the study of polymers, random polygonal paths in $\R^3$ also play a prominent role. Again we create a closed self-avoiding path by joining $n$ straight segments, but these are now distributed according to some continuous law. Two common choices for the distribution of the segments are the \emph{equilateral} with uniform distribution on the 2-sphere, and the \emph{Gaussian} with standard $3$-normal distribution. 

\begin{figure}[H]
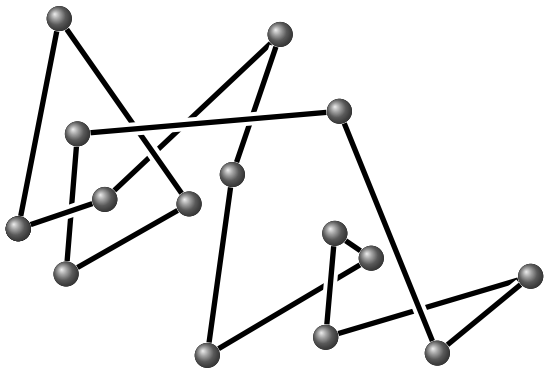

\begin{center}
\knotpoly
\caption{A closed polygon that realizes the trefoil in $\R^3$.}
\label{polygonal}
\end{center}
\end{figure}

That the walk is self-avoiding is usually satisfied with probability one, but more care is needed to make sure that the walk is closed. Details of this vary with the specific model and sampling method. We remain brief and only mention that it is possible to guarantee this in the Gaussian model by adding a constant drift.

It was conjectured by Frisch and Wasserman~\cite{frisch1961chemical} that polygonal walks are also unknotted with vanishing probability. Numerical simulations suggested exponential decay in $n$ for various different polygonal models~\cite{des1979topological, le1980monte, michels1982probability, michels1986topology, koniaris1991knottedness}. Diao, Pippenger, and Sumners~\cite{diao1994random} proved $\exp(-n^\varepsilon)$ for some $\varepsilon>0$ for Gaussian-steps polygons. This was extended to equilateral polygons~\cite{diao1995knotting} and other models~\cite{van2007knotting}.

General polygonal walks have an advantage over grid walks, in being space-isotropic. This is more realistic for polymers, and more robust to variations. Here to simulate effects of excluded volume constraints, one often replaces segments with rods and points with beads of positive radius. It is also interesting to consider polygons packed in a confined space such as a cube or a tube. Other variations of the model allow simulating bending rigidity, tension, pressure, thermodynamic entropy, and interaction between particles. See~\cite{micheletti2011polymers} for a thorough review.

For numerical experiments, such models are often approximately sampled via Markov chains in the configuration space, with a variety of local and global moves based on re-ordering, rotation, reflection, and more~\cite{alvarado2011generation}. There is currently much activity in search of faster rigorous sampling algorithms, with new techniques from symplectic geometry~\cite{cantarella2016symplectic,cantarella2016fast} and convexity~\cite{chapman2016ergodic}.

The resulting knots were classified for large samples in the various experiments. It turns out that, in several polygonal and grid models, the frequency at which a knot $K$ occurs is well approximated by $P[K_n{=}K] = C_K \left(n/N\right)^{\alpha_K} e^{-(n/N)}$. The constant $N$ depends only on the model, while for every knot $K$ the exponents $\alpha_K$ seem to be universal among different models~\cite{deguchi1994statistical,deguchi1997universality,orlandini1998asymptotics,millett2005universal,van2011universality}. Further experiments indicated that the $m$th most frequent knot appears with probability of order~$m^{-1.75}$~\cite{cantarella2016fast}.

\paragraph{Smoothed Brownian Motion} \label{smooth-model}

A substantially less studied subject is knotting from non-piecewise-linear three-dimensional random walks. A random polygon in the Gaussian model can be viewed as a linear interpolation between a finite number of points from a continuous Brownian bridge taken at constant time intervals. However, Brownian motion cannot model a knot as it is self-intersecting with probability one. Moreover, Kendall~\cite{kendall1979knotting} showed that it would contains infinitely many knots of all types, in the sense of being contained in such knotted tubes. Is there a smooth model that avoids these problems but captures the behavior of Brownian motion other than in small scale?

The \emph{worm-like loop}~\cite{grosberg2000critical} from polymer physics is a conituum model that takes curvature into account. A smooth closed curve in~$\R^3$ is given weight proportional to $\exp\left(-\ell \int \|\ddot{\mathbf{r}}\|^2 ds \right)$, where $\mathbf{r}(s)$ is its arc-length parametrization and $\ell$ is a typical length of persistence to bending. A more general model of statistical mechanics, designated for ribbons, takes care of the bending direction and persistence to twisting as well~\cite{kessler2003effect}. It is known how to approximately sample from this model for open paths but not for closed ones.

In the search of a more numerically accessible model for worm-like loops, Rappaport, Rabin and Grosberg~\cite{rappaport2006worm, rappaport2007differential} suggested the following mathematical model. One way to construct a Brownian bridge in $\R^3$ is by the following Fourier series, with $w_k=1$. $$ \mathbf{r}(t) = \sum_{k=1}^{\infty} \frac{w_k}{k} \left( \mathbf{Z}_k \cos kt + \mathbf{Z}_k' \sin kt \right) \;\;\;\;\;\;\;\; \mathbf{Z}_k, \mathbf{Z}_k' \sim \text{3-normal iid.} $$ To obtain a smooth approximation, one can truncate the sum by $w_k = 1_{k \leq n}$ or by $w_k = e^{-k/n}$. Computer simulations of the second choice show an exponential decay of the unknotting probability~\cite{rappaport2006worm}. It is interesting to observe that the cut-off factor $w_k = e^{-(k/n)^2}$ is equivalent to smoothing the Brownian motion by convolution with a narrow Gaussian, which seems to be an appealing choice. 

Recent works~\cite{westenberger2016knots,rivin2016random} study the case of polynomially decaying coefficients $w_k=k^{-\alpha}$, where $\alpha \in \R$. For $\alpha > 0.5$ they derive bounds on the expected crossing number of a random knot, and on the variance of the linking number of a random link.

The parametrization of knots by a finite sum of cosines yields~\emph{Fourier knots}~\cite{buck1994random, trautwein1995harmonic, kauffman1997fourier}. As shown by Lamm~\cite{lamm2012fourier}, every knot can be obtained by taking $x(t) = \cos( k_x t + \phi_x)$, $y(t) = \cos( k_y t + \phi_y)$, and a finite sum of such cosines for~$z(t)$. This was recently improved by Soret and Ville~\cite{soret2016lissajous}, who showed that a sum of two cosines is sufficient. Taking a single cosine, $z(t) = \cos( k_z t + \phi_z)$ defines the well-studied  \emph{Lissajous knots}~\cite{bogle1994lissajous, jones1998lissajous, lamm1997there, hoste2006lissajous}. In~\cite{boocher2009sampling} and \cite{rivin2016random} experiments on random Fourier and Lissajous knots are reported.

\paragraph{Random Jump} \label{jump-model}

In the above random walk models the typical step length is small compared with the diameter of the whole embedded path. Millet~\cite{millett2000monte} suggests polygonal models where each point $X_1,\dots,X_n \in \R^3$ is independently sampled from some distribution, such as the uniform distribution on the cube $[0,1]^3$, or a spherically symmetric distribution with a uniform radius in $[0,1]$. To this end any rich enough distribution that almost surely avoids self-intersections will do, such as the 3-normal distribution, or uniform on the unit sphere~\cite{54412}. 

\begin{figure}[bt]
\begin{center}
\renewcommand{\arraystretch}{2}
\begin{tabular}{ccc}
\jumpknot && \jumplink \\
A knot in the unit ball && A two-component link in the cube
\end{tabular}
\caption{The random jump model.}
\label{figjumps}
\end{center}
\end{figure}

By sampling $X_1,\dots,X_n$ and $Y_1,\dots,Y_m$ independently with the same 3-dimensional distribution, the above extends to two-component links~\cite{arsuaga2007linking}, and similarly for any number of components.

We still do not know how likely it is to encounter the unknot in the random jump model. Numerical experiments indicate that this probability vanishes faster than $\exp(-O(n))$~\cite{millett2000monte, arsuaga2007sampling}. This provides evidence for a strong form of the above-mentioned Delbruck--Frisch--Wasserman conjecture in this model. Similar conclusions seem to apply to any fixed knot. Experiments with the cube model suggest that the expected knot \emph{determinant} is $\omega\left(\exp n^2\right)$. It is proposed in~\cite{arsuaga2007sampling} that most knots in this model are prime. It was suggested~\cite{54417} that the expected crossing number in the spherical case is $\Theta(n^2)$.

Consider the \emph{linking number} $L_{mn}$ of a random two-component link with $n$ and $m$ segments. It is known~\cite{arsuaga2007linking,flapan2016linking} that its variance is $\Theta({nm})$, and it is conjectured that $L_{mn}/\sqrt{nm}$ converges in distribution to a Gaussian~\cite{panagiotou2010linking,karadayi2010topics}. Based on our analysis of the Petaluma model~\cite{even2016invariants} we tend to doubt this conjecture. Rather, we suspect that the tails of the limit distribution decay exponentially.

A more symmetric variant of the random jump model has been suggested~\cite{daniwise}, which takes place in $S^3$ visualized as the unit sphere in $\R^4$. A sequence of uniformly random points can be connected along the geodesics, which are the great circles.

\medskip

These \emph{random jump} or \emph{uniform random polygon} (URP) models, were originally proposed to illustrate the effect of spatial constraints on knotted molecules~\cite{millett2000monte}. In some bacteriophages, for example, a circular DNA molecule is densely packed inside a spherical capsid. Experiments show that more complex knots are likely to be produced, compared to unconstrained DNA of similar length in free solution~\cite{arsuaga2002knotting}. The observed distribution is also biased towards chiral knots and especially torus knots~\cite{arsuaga2005dna}. 

The explanation of these findings requires more realistic simulations that take into account various biophysical features, see e.g.~\cite{micheletti2008simulations, marenduzzo2009dna}. However, the simplicity of the random jump model makes it amenable for rigorous mathematical analysis, while it is arguably a prototype of a polygonal model in spatial confinement~\cite{arsuaga2007sampling}.

\paragraph{The Petaluma Model} \label{petaluma-model}

We now abandon random polygons, and move to more combinatorially oriented models. We start with the Petaluma model, studied by the author and collaborators~\cite{even2016invariants,even2017writhe,even2017petaluma}.

\begin{figure}[thb]
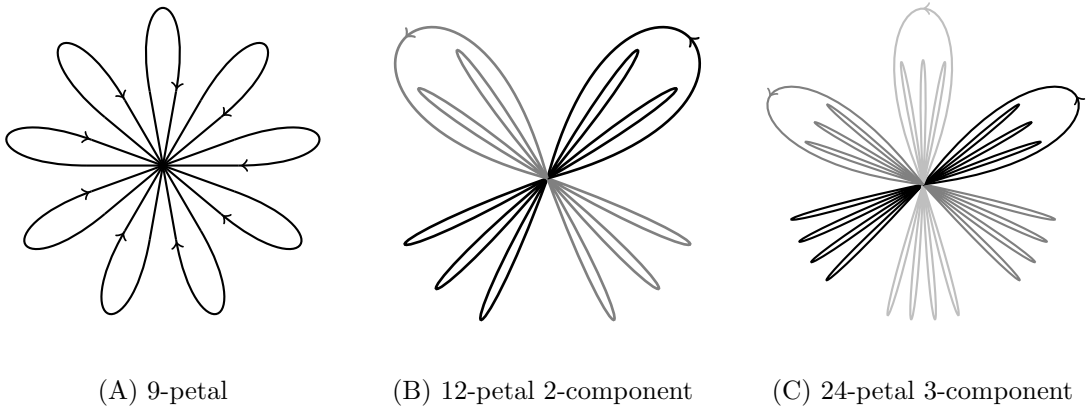

\begin{center}
\renewcommand{\arraystretch}{2}
\begin{tabular}{ccc}
\knotpetal & \;\;\linkpetal\;\; & \trilinkpetal \\
(A) $9$-petal & (B) $12$-petal $2$-component &  (C) $24$-petal $3$-component
\end{tabular}
\caption{Petal diagrams for knots and links.}
\label{petaldiagrams}
\end{center}
\end{figure}

Adams et al.~\cite{adams2015knot,adams2015bounds} have shown that every knot or link can be positioned so that its planar projection is injective except for a single point. Several projected strands may smoothly traverse this point of the \emph{\"{u}ber-crossing projection}, each originating at a different height. Moreover, every knot has a \emph{petal projection}, where the loops that emanate from the multi-crossing point have disjoint interiors. Consequently, petal projections are represented by a rose-shaped curve with an odd number of petals, as in Figure~\ref{petaldiagrams}A. 

In order to reconstruct the original knot we need only the relative ordering of the heights of the strands above the multi-crossing point. This information can be encoded by a permutation $\sigma \in S_{2n+1}$. We generate a random knot $K_{2n+1}$ in the \emph{Petaluma} model by picking $\sigma$ uniformly at random~\cite{even2016invariants}. By the construction of Adams et al., every knot $K$ is obtained with positive probability for $n$ large enough.

The Petaluma model extends to $k$-component links, by considering petal diagrams with $k$ components as in Figure~\ref{petaldiagrams}. In~\cite{even2017writhe} we study its extension to framed knots, which can be thought as knotted oriented ribbons.

In~\cite{even2016invariants,even2017writhe} we explicitly find the limiting distribution of the linking number of a two-component link, as well as the limiting distribution of the writhe of a random framed knot. We similarly present formulas for the moments of the Casson invariant $c_2$ and another finite type invariant appearing in the Jones polynomial. We elaborate on finite type invariants of random knots in the Petaluma model in Section~\ref{new-section} below.

As we show in a recent paper~\cite{even2017petaluma}, every particular knot appears in this model with vanishing probability. We conjecture that this probability decays at least exponentially with $n$, but currently the best bounds we have are $\Omega(n^{-n}) \leq P[K_{2n+1}{=}K] \leq O(n^{-0.1})$. 

It is of interest to understand the relation between the crossing number $c(K)$ of a knot and the least number of petals $p(K)$ needed to represent it. We show in~\cite{even2017petaluma} that $p(K) \leq O(c(K))$, and this bound is tight by results of Adams et al.~\cite{adams2015knot}. They have also shown that $c(K) \leq O(p^2(K))$, which is also tight.

Numerical simulations for $n \leq 100$ suggest that most knots in the Petaluma model are prime, and even hyperbolic. See Section~\ref{vol-exp} for more details, and further results by Adams and Kehne \cite{adams2017bipyramids,adams2016bipyramid,uberluma}. They went on to extend the Petaluma model to the \emph{\"{U}berluma} which contains all diagrams of one multi-crossing, allowing for nested loops.

\paragraph{Random Grid Diagrams} \label{grid-model}

Grid diagrams are a useful kind of regular knot diagrams. They describe all knots and links in a simple way~\cite{brunn1897uber,cromwell1998arc}. A \emph{grid diagram} consists of $n$ horizontal segments and $n$ vertical segments, where vertical segments always pass over horizontal ones. Each of the integers in $\{1,\dots,n\}$ appears as the $x$-coordinates of exactly one vertical segment. Likewise for the $y$-coordinates of the horizontal segments. 

\begin{figure}[bth]
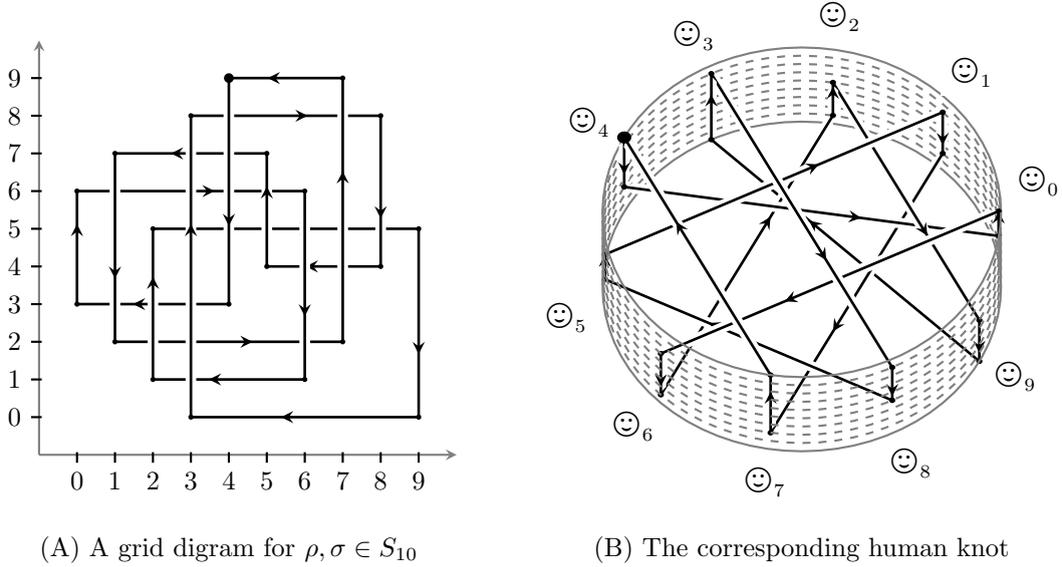

\begin{center}
\renewcommand{\arraystretch}{2}
\begin{tabular}{cccc}
\knotgrid &&& \humanknot \\
(A) A grid digram for $\rho,\sigma \in S_{10}$ &&& (B) The corresponding human knot
\end{tabular}
\caption{Here $\rho=(4,0,6,2,9,3,8,5,1,7)$ and $\sigma=(9,3,6,1,5,0,8,4,7,2)$.}
\label{griddiagram}
\end{center}
\end{figure}

A grid diagram is encoded by a pair of permutations $\rho,\sigma \in S_n$ for these horizontal and vertical coordinates respectively. We alternately take steps of the form $(\rho_{i},\sigma_{i}) \to (\rho_{i},\sigma_{i+1}) \to (\rho_{i+1},\sigma_{i+1})$ and so on. See~\cite{even2016invariants} for more details, and Figure~\ref{griddiagram}A for an example.

A random knot in the \emph{random grid} model is obtained by taking $\rho$ and $\sigma$ independently uniformly at random. Extensions to $k$-component links are easy and we omit further details. A~similar model that produces links of varying number of components was considered in a scheme for quantum money~\cite{farhi2012quantum}. 

We numerically compare the distribution of $c_2$ for the Petaluma and grid models, and find that they share many features, see Section~\ref{discussion}. As observed in~\cite{adams2015knot}, the Petaluma model is contained in the grid model, and obtained by conditioning on $\rho(k)=nk\bmod(2n+1)$.

Some preliminary work on precise moments' computation for finite type invariants in the random grid model has been done by Gal Lavi, Tahl Nowik, and the author~\cite{gallavi}. We report that $E[c_2] = n^2/288 + O(n)$ and $V[c_2] = 7n^4/194400 + O(n^3)$, which are of the same orders as in the Petaluma model, cf.~Section~\ref{new-section}. 

Two grid diagrams of the same knot can be related by a finite sequence of \emph{Cromwell moves}, which are local operations of three types, similar to the Reidemeister moves~\cite{cromwell1995embedding}. Witte et~al.~\cite{witte2016randomly} estimate the average writhe of a knot over its $n \times n$ grids, using a Markov chain of these moves. See also~\cite{farhi2012quantum}.

\medskip

We find a nice interpretation of the grid model in a common group-dynamic game named \emph{the human knot}~\cite{adams1994knot}. A group of $n$ two-handed participants stand in a circle. Each player chooses the next one at random and then they hold hands, until the last player holds the free hand of the first one. Their goal is to simplify the knot to a circle without letting their hands go, which is of course not always possible. 

To analyze this game, we introduce the assumption of transitivity. Namely, connected pairs of hands are ordered from bottom to top. See Figure~\ref{griddiagram}B, where the players correspond to axial segments on a cylinder, and connections are horizontal chords at different heights. If this ordering is uniformly random, then this construction is equivalent to a random grid diagram. Horizontal and vertical segments correspond to chords and players respectively. The permutation $\rho$ records the order at which players are connected, and $\sigma$ represents the relative order of the hands' heights. 

A related model, based on the human knot game, was suggested by Gilad Cohen~\cite{giladcohen}, who conducted computer experiments to study the distribution of the resulting knots.

\paragraph{Random Planar Diagrams} \label{diagram-model}

Planar diagrams are routinely used to represent knots and to investigate them. Naturally, this suggests the study of random knots by sampling diagrams with a given number of crossings. Such models were studied by several authors~\cite{schaeffer2004asymptotic,diao2005generating,diao2010mean,dunfield2014random,cantarella2016knot}, with various sampling methods.

\begin{figure}[bht]
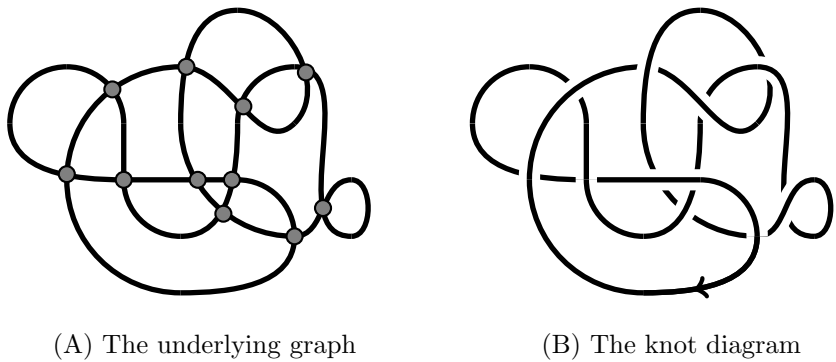

\begin{center}
\begin{tabular}{ccc}
\curvegraph && \curveknot\\
(A) The underlying graph && (B) The knot diagram
\end{tabular}
\caption{A random assignment of crossings to an $11$-vertex $4$-regular plane graph.}
\label{planar}
\end{center}
\end{figure}

To this end, we start with a generic smooth immersion of $S^1$ into $\R^2$ with $n$ traverse double points, considered up to diffeomorphism of the plane, as in Figure~\ref{planar}. This yields a $4$-regular plane graph, where loops and multiple edges are allowed.  Then each vertex is assigned either of the two possible crossing signs.

The number of $n$-vertex $4$-valent graphs in $\R^2$ is asymptotically exponential in $n$. However, an algorithm by Schaeffer~\cite{schaeffer2004asymptotic, brinkmann2007fast} uniformly samples such graphs with a base point, by generating a random rooted binary tree and matching leaves to non-leaves in some clever way. Some of the resulting graphs correspond to curves with several components, which is a problem if one is interested only in knots rather than links. One can either reject~\cite{dunfield2014random,cantarella2016knot} these curves, or modify~\cite{diao2005generating,diao2010mean} them, but this, however, ruins uniformity.

Some delicate issues of symmetry arise. Namely, do we care about orientation and mirror images? Should we distinguish between different planar diagrams which are equivalent in the sphere $S^2$? Do we want a base point on some edge? Finally, are different $n$-vertex graphs to be weighted equally or according to the number of non-equivalent diagrams they give rise to, which might be smaller than $2^n$ due to symmetries? However, all subtleties of this sort become negligible as $n$ grows~\cite{richmond1995almost,chapman2016asymptotic2}.

A recent advance in the study of this model is the establishment of a pattern theorem for diagrams by Harrison Chapman~\cite{chapman2016asymptotic,chapman2016asymptotic2}. This extends pattern theorems for planar maps~\cite{bender1992submaps}, and parallels the above-mentioned results for grid and polygonal knots. Chapman showed that small sub-diagrams appear $\Theta(n)$ times in an $n$-crossing knot or link diagram, except for an exponentially small probability. In particular, as $n$ grows the diagram contains a $3$-crossing trefoil summand and is hence nontrivial with high probability. Similar results hold if one restricts to \emph{prime} diagram, ones whose underlying graph is $4$-edges-connected.

Numerical experiments tell us more. Dunfield, Obeidin et al.~\cite{dunfield2014random,obeidin2016volumes} study random links, knots, and prime connected summands of knots in this model.  Their results suggest that several invariants, including the hyperbolic volume, grow linearly with $n$. Cantarella, Chapman and Mastin~\cite{cantarella2016knot,chapman2016asymptotic2} precisely compute knot probabilities for $n \leq 10$, and study their behavior for larger $n$ based on random samples. The methods used in these experiments are implemented into publicly available software packages: \emph{plCurve}~\cite{plcurve} and \emph{SnapPy}~\cite{snappy}.

\paragraph{Random Planar Curves} \label{curves-model}

Other models generate a random $4$-regular plane graph in various ways, and then assign crossing signs uniformly at random. For example, Diao et al.~\cite{diao2010mean} randomly add $n$ non-intersecting chords inside and outside an $n$-vertex cycle, to make it $4$-regular, and then toss a coin to decide each crossing. 

In the following random-crossing constructions the underlying graph is generated by sampling polygonal curves in the plane.

\begin{itemize}
\item Equilateral closed polygons in $\R^2$~\cite{michels1989distribution}.
\item Closed SAW in $\Z^2$ with diagonal crossings: \squarecross{(0)--(3)}{(1)--(2)}\; or \squarecross{(1)--(2)}{(0)--(3)}\;~\cite{guitter1999monte}.
\item Jumps between uniform points in the square $[0,1]^2$~\cite{arsuaga2007sampling,diao2010mean}.
\item A chain of chords between uniform points around the circle~\cite{giladcohen}.
\item The \emph{griddle}: Random grid diagrams with randomized crossings~\cite{workinprogress}.
\end{itemize}

There are close connections between the finite type invariants of such knots and those of the underlying curve~\cite{polyak1998invariants}. For example, the expected value of the Casson invariant $c_2$ is one eighth the \emph{defect}, a first-order invariant of the curve. In the griddle model we calculated $E[c_2] = E[\text{defect}]/8 = n^2/144 + O(n)$ and $V[c_2] = n^4/7776 + O(n^3)$, though $V[\text{defect}] = 29n^3/4050 + O(n^2)$~\cite{workinprogress}.

Finally, we note that given a $4$-valent graph in the plane, exactly two sign assignments produce an \emph{alternating link diagram}, where over-crossings and under-crossings alternate as one travels along the link. Diao et al.~\cite{diao2005generating,arsuaga2007sampling,diao2010mean} and Obeidin~\cite{obeidin2016volumes} used this observation to construct models for prime alternating knots and links. Except for the $(2,n)$-torus these are hyperbolic links, whose volume can be read off the diagram up to a multiplicative constant~\cite{lackenby2004volume}. Taking the uniform distribution over prime alternating link diagrams, the expected hyperbolic volume is linear in the crossing number~\cite{obeidin2016volumes}.

\paragraph{The Knot Table Model} \label{table-model}

\begin{figure}[tbh] 
\center
\includegraphics[trim={1.75cm -1.5cm 1.75cm -1cm},clip,width=\columnwidth]{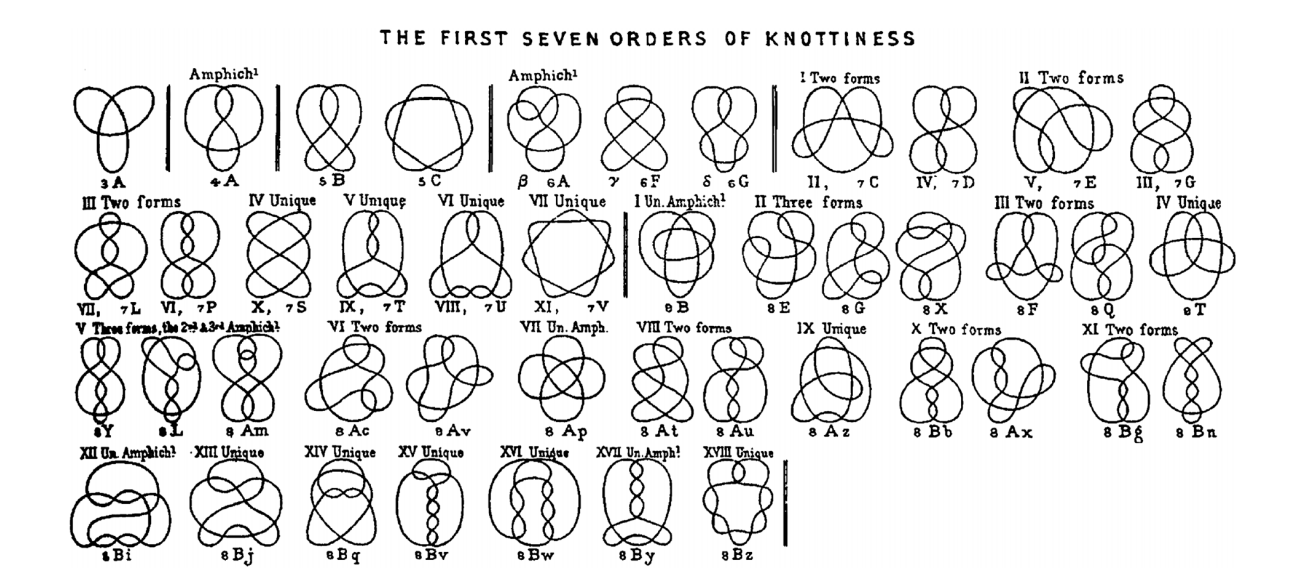}
\captionsetup{width=0.8\textwidth}
\caption{Excerpt from Tait's original table of knots with up to $8$ crossings~\cite{tait1884first}. Note that unlike the discussed model it contains only alternating knots, with several equivalent diagrams for some of them. \\ ~}
\label{tait}
\end{figure}

The crossing number is perhaps the most popular measure for knot complexity. Historically, prime knots are tabulated and nomenclated according to their crossing number, as reflected in the widely used Alexander–-Briggs--Rolfsen knot notation~\cite{alexander1926types, rolfsen1976knots}. See also Figure~\ref{tait}.

Consequently, many investigators find it quite natural to generate random prime knots by uniformly sampling from knot tables with up to $n$ crossings. If one cares about chirality and orientation, these can be decided by further coin flips.

It is known that there are exponentially many knots with $n$ crossings \cite{ernst1987growth,welsh1991number,sundberg1998rate}, but the exact count is known only for small $n$~\cite{hoste1998first}. The difficulties in recognition and enumeration of $n$-crossing knots make this model less suitable for precise computations, though it is known that most knots are not rational~\cite{ernst1987growth}, nor are most links alternating~\cite{thistlethwaite1998structure}. 

The vast majority of knots with up to $n \leq 16$ are hyperbolic, which may suggest that their asymptotic proportion tends to $1$. This is however not likely to be true, in view of a recent surprising result of Malyutin~\cite{malyutin2016question}. He assumes the plausible, but still unproven, conjecture that the crossing number is weakly monotone with respect to connected sum. The crux of his proof is the addition of small satellite configurations to existing diagrams.

\paragraph{Random Braids} \label{braid-model}

\begin{figure}[H]
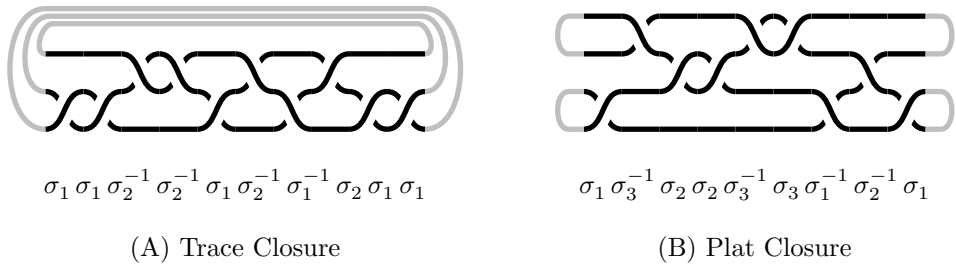

\begin{center}
\renewcommand{\arraystretch}{2}
\begin{tabular}{ccc} 
\braidknot && \bridgeknot \\
$\sigma_1^{}\,\sigma_1^{}\,\sigma_2^{-1}\,\sigma_2^{-1}\,\sigma_1^{}\,\sigma_2^{-1}\,\sigma_1^{-1}\,\sigma_2^{}\,\sigma_1^{}\,\sigma_1^{}$ && 
$\sigma_1^{}\,\sigma_3^{-1}\,\sigma_2^{}\,\sigma_2^{}\,\sigma_3^{-1}\,\sigma_3^{}\,\sigma_1^{-1}\,\sigma_2^{-1}\,\sigma_1^{}$ \\
(A) Trace Closure && (B) Plat Closure 
\end{tabular}
\caption{Random knots in the braid model.}
\label{braid}
\end{center}
\end{figure}

It goes back to Alexander that every knot or link is the closure of some \emph{braid}~\cite{lickorish1997introduction}. Namely, it can be represented by some $m$ intertwining strings that monotonously go from left to right, and close at some canonical way as in Figure~\ref{braid}. Such braids form a group~$B_m$, with generators $\{\sigma_i^{\pm 1}\}_{1 \leq i < m}$ that correspond to swapping strings $i$ and $i+1$, and appropriate relations. 

There is recent interest in generating knots by random walk in the braid group. This parallels well-known constructions of random 3-manifolds~\cite{dunfield2006finite} and more.

Such a model is defined in terms of a probability distribution on a finite subset of the braid group $B_m$, such as the generators $\sigma_i^{\pm 1}$. A random knot is obtained by $n$-step random walk in these generators, with some standard closure as depicted in Figure~\ref{braid}. The context of Markov Chains on groups proves useful in the analysis of this model~\cite{nechaev1996random}. 

This definition yields random links of a varying number of components. For fixed $m$ and large $n$ we obtain knots with probability about $1/m$. Additionally, only links of \emph{braid index} or \emph{bridge index} at most $m$ appear, according to the closure convention. Remarkably, random knots and links in this setting are hyperbolic with high probability~\cite{malyutin2012quasimorphisms,ma2013components, ma2014closure, ito2015structure, ichihara2015most, ichihara2016random}.

\paragraph{Crisscross Constructions} \label{cross-model}

This family of random models includes several constructions in which a planar curve is explicitly specified, and all randomness comes from the choice of crossing signs, sampled independently and uniformly at random. 

\medskip

\begin{figure}[H]
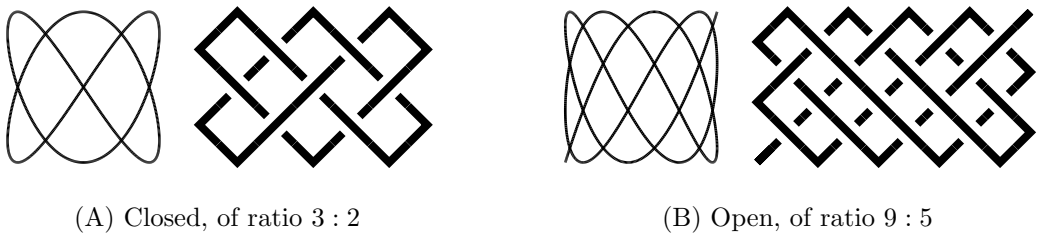

\begin{center}
\renewcommand{\arraystretch}{2}
\begin{tabular}{ccc}
\lissajouscurve \;\;\; \lissajousknot & \;\;\;\;\;\;\;\;\;\; & \harmoniccurve \;\;\; \harmonicknot \\
(A) Closed, of ratio $3:2$  & &  (B) Open, of ratio $9:5$
\end{tabular}
\caption{Billiard table diagrams from Lissajous curves.}
\label{billiard}
\end{center}
\end{figure}

One source for such models is planar Lissajous curves~\cite{lissajous1857memoire}, illustrated in Figure~\ref{billiard}. These closed curves are parametrized by $(\cos(at + \phi)$, $\cos bt)$ where $t \in [0, 2\pi]$ with ratio $b\,{:}\,a \in \Q$ and a phase shift $\phi \in \R$. We also consider the \emph{open}  curve $(\cos at, \cos bt)$ where $t \in [0, \pi]$, being closed from the outside. These curves are plane isotopic to the polygonal trajectory of a billiard ball in~$[0,1]^2$, fired at slope~$b/a$~\cite{jones1998lissajous}.

The three-dimensional analogues of these curves constitute \emph{Lissajous knots}~\cite{bogle1994lissajous, lamm1997there, jones1998lissajous} and \emph{Harmonic Knots}~\cite{comstock1897real, koseleff2011chebyshev}, but these families do not contain all knots. However, planar Lissajous curves with suitable crossing signs do give rise to all knots. This underlies the construction of the above-mentioned \emph{Fourier Knots}~\cite{buck1994random, trautwein1995harmonic, kauffman1997fourier, hoste2006lissajous, lamm2012fourier, soret2016lissajous} and \emph{Chebyshev Knots}~\cite{koseleff2011chebyshev}, as well as the next random construction, the \emph{billiard table model} suggested by Cohen and Krishnan~\cite{cohen2015random}.

A random knot $K_{b:a}$ is thus obtained by randomizing the crossing signs, as in Figure~\ref{billiard}. It can also be regarded as a special case of the random braid model. For example, the case $a=5$ as in Figure~\ref{billiard}B is generated by the $16$ elements $\{\sigma_1^{\pm}\sigma_3^{\pm}\sigma_2^{\pm}\sigma_4^{\pm}\}$ with the uniform distribution.

In~\cite{cohen2016crossing} we study the asymptotic properties of $K_{n:3}$, which yields random two-bridge knots, also known as \emph{rational knots}~\cite{kauffman2004classification}. We show that the probability of obtaining any particular knot is $(\alpha+o(1))^n$ for $\alpha = \sqrt[3]{27/32} \approx 0.945$, and the crossing number is $(\beta+o(1))n$ in probability, for $\beta = (\sqrt{5}-1)/4 \approx 0.309$. 

We remark that, without restricting to fixed diagrams, other random models arise from the highly developed theory of rational knots. In particular, a random braid in $\{\sigma_1,\sigma_2^{-1}\}^{\star} \subset B_4$ yields a rational knot by its Conway symbol~\cite{conway1970enumeration}. See~\cite{ernst1987growth} and \cite{diao2010mean} for corresponding results.

\medskip

\begin{figure}[t]
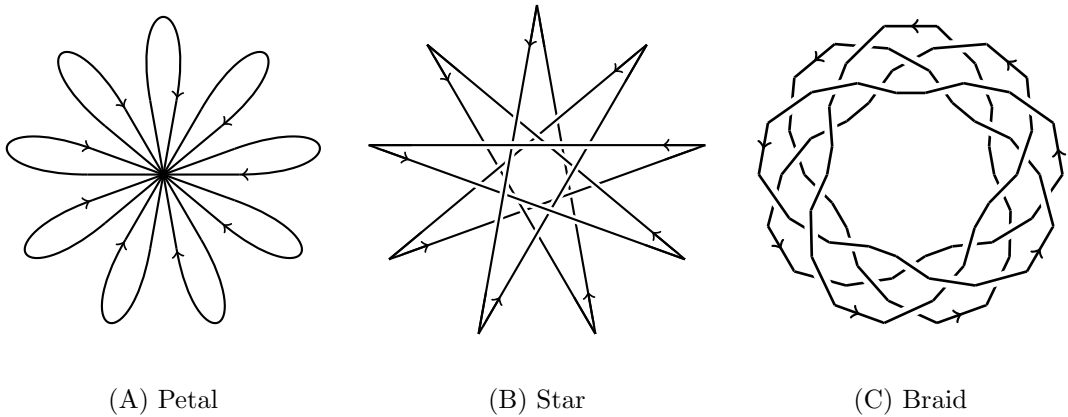

\begin{center}
\renewcommand{\arraystretch}{2}
\begin{tabular}{ccc} 
\knotpetal & \starknot & \starbraid \\
(A) Petal & (B) Star & (C) Braid
\end{tabular}
\caption{From petal diagrams to regular knot diagrams.}
\label{petalstarbraid}
\end{center}
\end{figure}

Star diagrams are obtained from $(2n+1)$-petal diagrams by straightening the segments between petal tips. See Figure~\ref{petalstarbraid}A-B. A random knot in the \emph{star model} is generated by randomizing the $(n-1)(2n+1)$ crossings. Star diagrams are plane isotopic to closed $n$-braids~\cite{adams2015knot}, as demonstrated in Figure~\ref{petalstarbraid}B-C.

The star model yields all knots, since the Petaluma model does, but with quite different distribution. We show in~\cite{even2016invariants} that its expected Casson invariant is $E[c_2] = n^3/12 + O(n^2)$ with a standard deviation of $n^2/\sqrt{24} + O(n^{3/2})$. This means that $c_2$ drifts away from zero. 

Chang and Erickson~\cite{chang2015electrical} consider a generalization of the star model. They define the \emph{flat torus} diagram $T(p,q)$ as the closed braid $(\sigma_1 \sigma_2 \cdots \sigma_{p-1})^q$, and assign crossing signs at random. The star model is $T(n,2n+1)$, as shown in Figure~\ref{petalstarbraid}C for $n=4$. Following Hayashi et al.~\cite{hayashi2012minimal}, they show that the expected Casson invariant of $T(n+1,n)$ is $\Theta(-n^3)$. It is conceivable that this latter model contains all knots as well. 

\medskip

The probability space in such crisscross models consists of $2^c$ crossing states. Some invariants are more accessible in this simple setting, as they are computable by summation over $2^c$ local configurations at the $c$ crossings. One important example is the Kauffman Bracket~\cite{kauffman1987state}, and its connections to statistical physics~\cite{kauffman1988statistical,jones1989knot,wu1992knot}. 

For crisscross diagrams on the 2-dimensional lattice, rather similar to the above ones, the degree distribution of the Jones polynomial is analyzed in terms of the Potts model from statistical mechanics~\cite{grosberg1992algebraic, nechaev1996statistics, vasilyev2001thermodynamics}.

\paragraph{Miscellanea}

We have attempted to cover the main themes of random knot models. Of course, our list of models and results is not completely exhaustive, neither historical, and to some extent reflects our own viewpoint. To conclude, we mention some random ideas in further directions.

\smallskip

Various models from the natural sciences seek to emulate dynamical processes of knot formation in real life scenarios. Some studies describe numerical simulations of a polygonal DNA chain that folds, coils and spools within a cavity, before its two ends anneal and produce a knot~\cite[for example]{arsuaga2008dna, marenduzzo2009dna}. Such dynamical models are important for understanding biological processes by comparing simulated and observed data, but usually they don't lend themselves easily to mathematical analysis.

\smallskip

Other studies~\cite{flammini2004simulations, hua2007random, liu2008efficient, szafron2011knotting, cheston2014new} are inspired by the interaction between DNA and \emph{topoisomerase}, a specific enzyme that cuts and rejoins strands, and thus modifies their topological state. Such strand-passage models induce transition probabilities between knot types, which can be estimated by numerical simulations, and these lead to a stationary equilibrium distribution over knots.

\smallskip

Finally, Babson and Westenberger study knots obtained from a curve in $\R^n$ by projecting to~$\R^3$ in a random direction. They relate several of the above constructions to this original framework~\cite{westenberger2016knots}. 

\smallskip

In principle, any reasonable way to construct or represent knots could be turned into a random model. Another case in point are trajectories of dynamical systems, such as three-dimensional billiard~\cite{jones1998lissajous}.

\medskip

\section{A Closer Look at the Petaluma Model}\label{new-section}

We now focus on random knots and links in the Petaluma model~(\ref{petaluma-model}), and discuss the distribution of their finite type invariants and hyperbolic volume. First we recall the definition of finite type invariants, given in terms of singular knots and links~\cite{birman1993knot}.

\paragraph{Finite type invariants}\label{fti-def}

Unlike a regular knot, which is a smooth embedding of~$S^1$ into~$\R^3$ up to isotopy, a \emph{singular} knot is allowed to have finitely many double points of transversal self intersection. Each of these points can be locally \emph{resolved} in two well-defined ways: positive~\tikz[thick,->,baseline=-2]{\draw[black](0.3,-0.1) -- (0,0.2);\draw[white,-,line width=3](0.03,-0.1) -- (0.33,0.2);\draw[black](0.03,-0.1) -- (0.33,0.2);}, and negative~\tikz[thick,->,baseline=-2]{\draw[black](0.03,-0.1) -- (0.33,0.2);\draw[white,-,line width=3](0.3,-0.1) -- (0,0.2);\draw[black](0.3,-0.1) -- (0,0.2);}.

Let $v$ be a knot invariant taking values in some abelian group, usually in $\Z$. The extension of $v$ to singular knots is given by $v(K) = v(K_p^+) - v(K_p^-)$, where $K_p^{\pm}$ are the two resolutions of the singular knot $K$ at the double point $p$. By recursion, the value of $v$ on a singular knot with $m$ double points is given by a signed sum of its value on $2^m$ regular knots. We say that $v$ is a \emph{finite type} knot invariant of \emph{order} $m$ if it vanishes on all singular knots with $m+1$ double points. 

This condition is satisfied by several well-studied knot invariants, such as coefficients of knot polynomials~\cite{bar1995polynomial,chmutov2012introduction} and the Kontsevich integral~\cite{bar1995vassiliev, chmutov2001kontsevich}. There is only one knot invariant of order two, up to affine equivalence -- the \emph{Casson invariant} $c_2(K)$, which is the coefficient of $x^2$ in the Alexander--Conway polynomial $C_K(x)$. It similarly appears in the modified Jones polynomial, $V_K(e^x)$ considered as a power series in $x$, which also yields an invariant $v_3(K)$ of order three. The number of new independent finite type invariants grows with the order: $3$ invariants of order four, $4$ of order five, $9$ of order six, etc.~\cite{bar1995vassiliev}. 

No invariant of knots has order one. However, the Gauss \emph{linking number} $lk(L)$ is a classical first order invariant of two-component links. Also the framing number, or \emph{writhe} $w(K)$ as in~\cite{even2017writhe}, is a first order invariant of framed knots.

\paragraph{Asymptotic Distributions}

~Finite type invariants of random knots and links in the Petaluma model (\ref{petaluma-model}) have been studied by Hass, Linial, Nowik, and the author~\cite{even2016invariants,even2017writhe}. In particular, we have investigated how these invariants scale and distribute for knots with a large number of petals. 

Consider the Casson invariant of a random knot with $2n+1$ petals. It is not hard to observe that $c_2(K_{2n+1}) = O(\pm n^4)$, which is shown to be sharp for torus knots and other explicit constructions. However, we have found that the typical order of magnitude of the Casson invariant is actually~$n^2$. Indeed, its expectation is $E[c_2] = n(n-1)/24$, its variance is $V[c_2] = 7/960 \cdot n^4 + O(n^3)$, and such formulas have been given for all moments, yielding $E[c_2^k] = \Theta(n^{2k})$. We find it intriguing that the distribution of the properly normalized Casson invariant $c_2/n^2$ is asymmetric and not centered at zero, asymptotically as $n \to \infty$. 

The third order invariant $v_3(K_{2n+1})$ is antisymmetric with respect to reflection, hence its distribution is symmetric around zero. As for its even-order moments, we have similarly shown $E[v_3^k] = O(n^{3k})$, e.g., $V[v_3] = 9298/5443200\cdot n^6 + O(n^5)$.

In terms of their moments, $c_2$ grows as $n^2$ and $v_3$ as $n^3$. This naturally suggests that an $m$th order invariant of random knots with $n$ petals asymptotically scales as~$n^m$. In~\cite{even2016invariants} we conjecture that $v_m(K_{2n+1})/n^m$ weakly converges to a limiting distribution as $n \to \infty$ for every finite type invariant $v_m$ of order $m$. The existence of continuous limit distributions for $c_2$ and $v_3$ is supported by computational evidence, as discussed below.

We have established such a limiting distribution in two cases: the linking number of a random two component link with $2n$ petals in each component, and the writhe of a random framed knot with $2n+1$ petals. Both are first order invariants, and obtain integer values sharply between~$\pm n^2$. In~\cite{even2016invariants} we prove that $lk(L_{2n,2n})/4n$ converges to the \emph{logistic distribution}, with density function $f(t)= \pi/\cosh^2(2\pi t)$. The normalized writhe $w(K_{2n+1})/n$ converges to another non Gaussian limiting distribution, established and described in~\cite{even2017writhe}. 

Our proofs combine the method of moments with careful combinatorial analysis of the limiting moments of these invariants, expressed via Gauss diagram formulas.

\paragraph{Numerical Experiments}\label{fish-section}

\begin{figure}[H] 
\center
\includegraphics[trim={0.5cm 0cm 3cm 0cm},clip,width=\columnwidth]{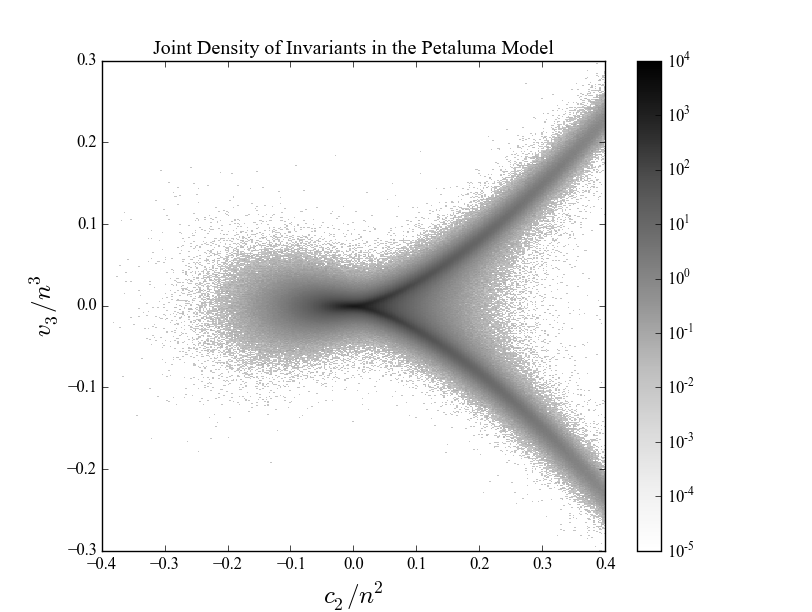}
\captionsetup{width=0.7\textwidth}
\caption{The normalized distribution of $c_2$ and $v_3$ for a random knot $K_{41}$, based on $10^8$~random samples. \\ ~}
\label{fish}
\end{figure}

We study the invariants $c_2(K_n)$ and $v_3(K_n)$ in the Petaluma model, by computing their values for a random sample of permutations in $S_n$. Comparing the results for various values of $n$, we observe that as $n$ grows the joint distribution of $c_2/n^2$ and $v_3/n^3$ seems to converge to a continuous bivariate distribution of a certain shape. The heat map in Figure~\ref{fish} shows the resulting density function of this distribution for $n=41$, which seems to be a good approximation of the conjectured limiting distribution.

\bigskip

\newpage	

The planar representation of these two invariants follows previous work by Willerton~\cite{willerton2002first,chmutov2012introduction} and Okuda~\cite{ohtsuki2002problems}, who generated scatter plots of $(c_2,v_3)$ for all prime knots with up to $n$ crossings. They similarly obtained fish-shaped figures, although it is unclear how these should scale as the crossing number grows. The Petaluma model may provide a more concrete way to catch this fish, in the form of a limit density function defined on $\R^2$.   

Besides representing the first two finite type invariants, the planar map $\varphi:K \mapsto (c_2(K),v_3(K))$ has some interesting properties. As observed by Dasbach et al.~\cite{dasbach2001quantum}, the evaluation of the Jones polynomial at roots of unity near $1$ can be approximated by $V_K(e^{ih}) = 1 + 3c_2h^2 + 6v_3h^3i + O(h^4)$, and this yields similar fish graphs for $V_K(e^{2\pi i/N})$ in the complex plane, for $N \gg n$. 

Note that by the multiplicativity of the Jones polynomial, the map $\varphi$ is additive with respect to connected sum: $\varphi(K\#K') = \varphi(K) + \varphi(K')$ in $\Z^2$. Using this fact and some known constructions one can show that as $n$ grows the resulting point set of all $(c_2/n^2,v_3/n^3)$ is dense in $\R^2$. We actually conjecture that the limiting bivariate distribution has positive density everywhere in the plane.

\paragraph{Hyperbolic Volume}\label{vol-exp}

We conclude this section with further numerical experiments, concerning the distribution of the hyperbolic volume in the Petaluma model, as approximated by the \emph{Sage} package \emph{SnapPy}~\cite{snappy}.

\bigskip

\begin{figure}[H] 
\center
\includegraphics[trim={1.5cm 0cm 1.5cm 0cm},clip,width=\columnwidth]{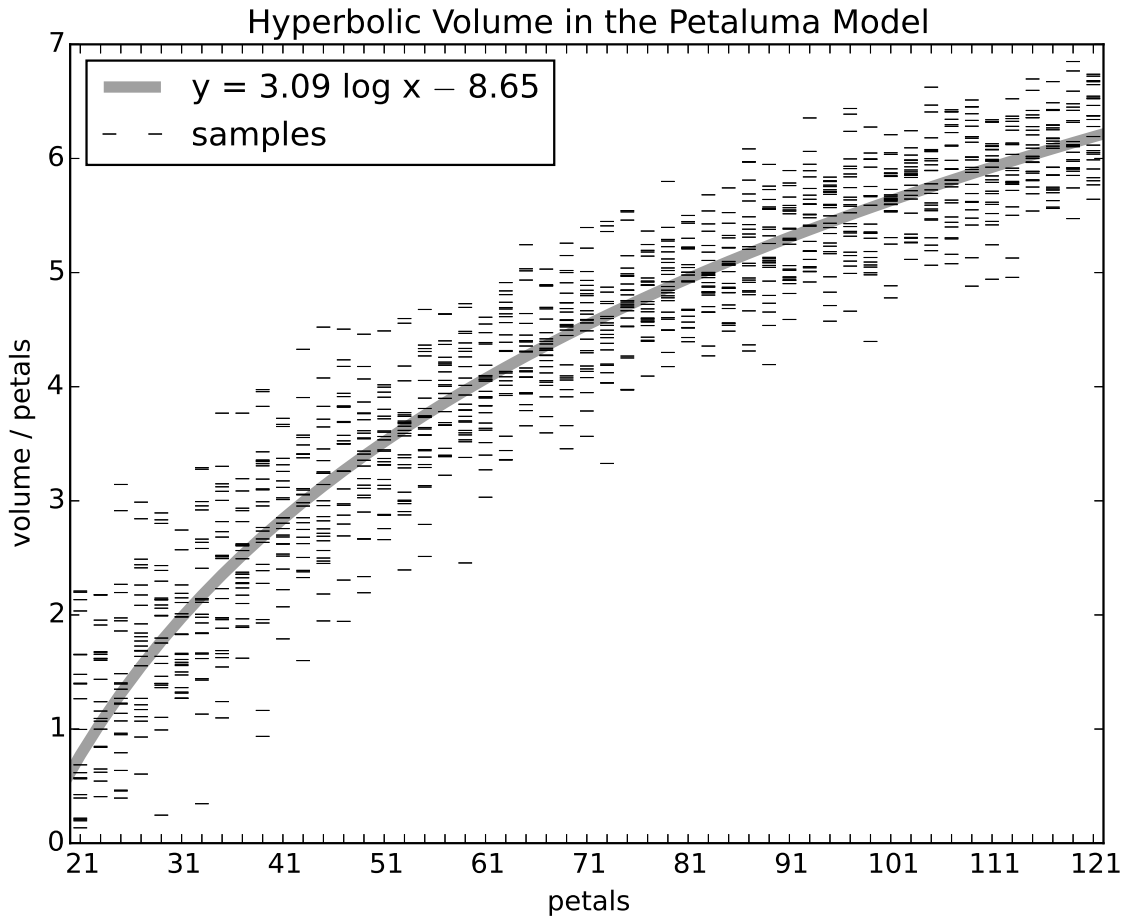}
\captionsetup{width=0.75\textwidth}
\caption{The hyperbolic volume per petal grows with the number of petals. This is based on random samples of $25$ knots with $21$ to $121$ petals. Non-hyperbolic knots ($<2.5\%$) were omitted. \\ ~}
\label{vol}
\end{figure}

\newpage

As mentioned in Section~\ref{petaluma-model}, our simulations show that randomly sampled knots with up to $200$ petals are mostly hyperbolic. This trend seems to strengthen with increasing number of petals, although one must be careful drawing conclusions from small cases, cf.~\cite{malyutin2016question} mentioned in Section~\ref{table-model}.

Figure~\ref{vol} shows how the empirical hyperbolic volume grows super-linearly with the number of petals. More speculatively, the volume of an $n$-petal knot appears to be concentrated around a curve of the form $An \log Bn$, which seemed to fit better than a linear function, or one of order $n^{3/2}$. These experiments have been repeated by Adams and Kehne~\cite{uberluma}. They have also proved that the expected volume is at most $4 \pi \, n \log n$, by constructing a pyramid decomposotion of the petal knot complement~\cite{adams2017bipyramids,adams2016bipyramid}. Any such lower bound would be of great interest.

\medskip

\section{Discussion}\label{discussion}

This great variety of approaches for random knot models suggests that we ask how they differ. Do they exhibit some kind of common properties? By what means should we compare models? What do they teach us about knot invariants and knot theory? Below we record some thoughts concerning these questions. 

\paragraph{Local knotting}\label{local-section}

The Delbruck--Frisch--Wasserman conjecture, that a typical random knot is non-trivial, has been proved by now in several models. Some insight on their properties can be gained by comparing the arguments involved in these proofs.

The knottedness of random polygonal and grid walks~(\ref{walk-model},\ref{polygonal-model}) is based on the fact that such knots tend to have many spatially \emph{localized} connected summands. This phenomenon can be attributed to the small steps taken in these models~\cite{sumners1988knots, pippenger1989knots, diao1994random, diao1995knotting}. We do know, however, that large scale knotting occurs as well~\cite{jungreis1994gaussian, diao2001global}. Also for planar diagrams~(\ref{diagram-model}), knottedness follows from the existence of small prime summands in random knot and link diagrams~\cite{chapman2016asymptotic2}. Even for prime knots in the knot table model~(\ref{table-model}), local configurations of a double figure-eight knot provide a satellite decomposition~\cite{malyutin2016question}.

In contrast to the highly composite knots produced by small-steps models, we believe that models of \emph{non-local} nature yield knots with much simpler factorization. By non-local we mean that the typical step length is comparable to the diameter of the whole curve. 

For example, local entanglements yield only a vanishing probability of order $1/n^3$ for a trefoil summand in the Petaluma model~(\ref{petaluma-model}). Indeed, its knottedness with high probability was shown by other means, a coupling argument based on the effect of random crossing changes on finite type invariants~\cite{even2017petaluma}. As mentioned above, numerical experiments indicate that these knots are mostly hyperbolic, so that any connected sum or satellite-type decomposition might become rare.

\paragraph{Dimension}\label{dimension-section}

It would be interesting to further distinguish knot models from each other by their asymptotic topological features. On the other hand, it would be very interesting to discover \emph{universal} phenomena and parameters that hold for a variety of different models.

We shall venture some speculations along these lines. As a first step, consider the following three classes of random models.

\begin{framed}
\begin{itemize}
\item[\textbf{1D}] Grid walks~(\ref{walk-model}), polygonal walks~(\ref{polygonal-model}), and smoothed Brownian motion~(\ref{smooth-model}).
\item[\textbf{2D}] Random planar diagrams~(\ref{diagram-model}), the griddle~(\ref{curves-model}), knot table~(\ref{table-model}), and star~(\ref{cross-model}).
\item[\textbf{3D}] Random jumps~(\ref{jump-model}), the Petaluma~(\ref{petaluma-model}), and grid diagrams~(\ref{grid-model}).
\end{itemize}
\end{framed}

This classification attempts to grasp the ``dimension'', or general shape, of the actual spatial curves constructed by the different models, in some loose and undefined sense. It is a fundamental challenge to characterize such a classification precisely.

Would it be possible to reconstruct the class to which some random model belongs, by looking only at the asymptotics of the topological invariants of the resulting knots?

\paragraph{Comparing Invariants}\label{casson-section}

\begin{figure}[H] 
\center
\includegraphics[trim={1.8cm 0cm 1.8cm 0cm},clip,width=\columnwidth]{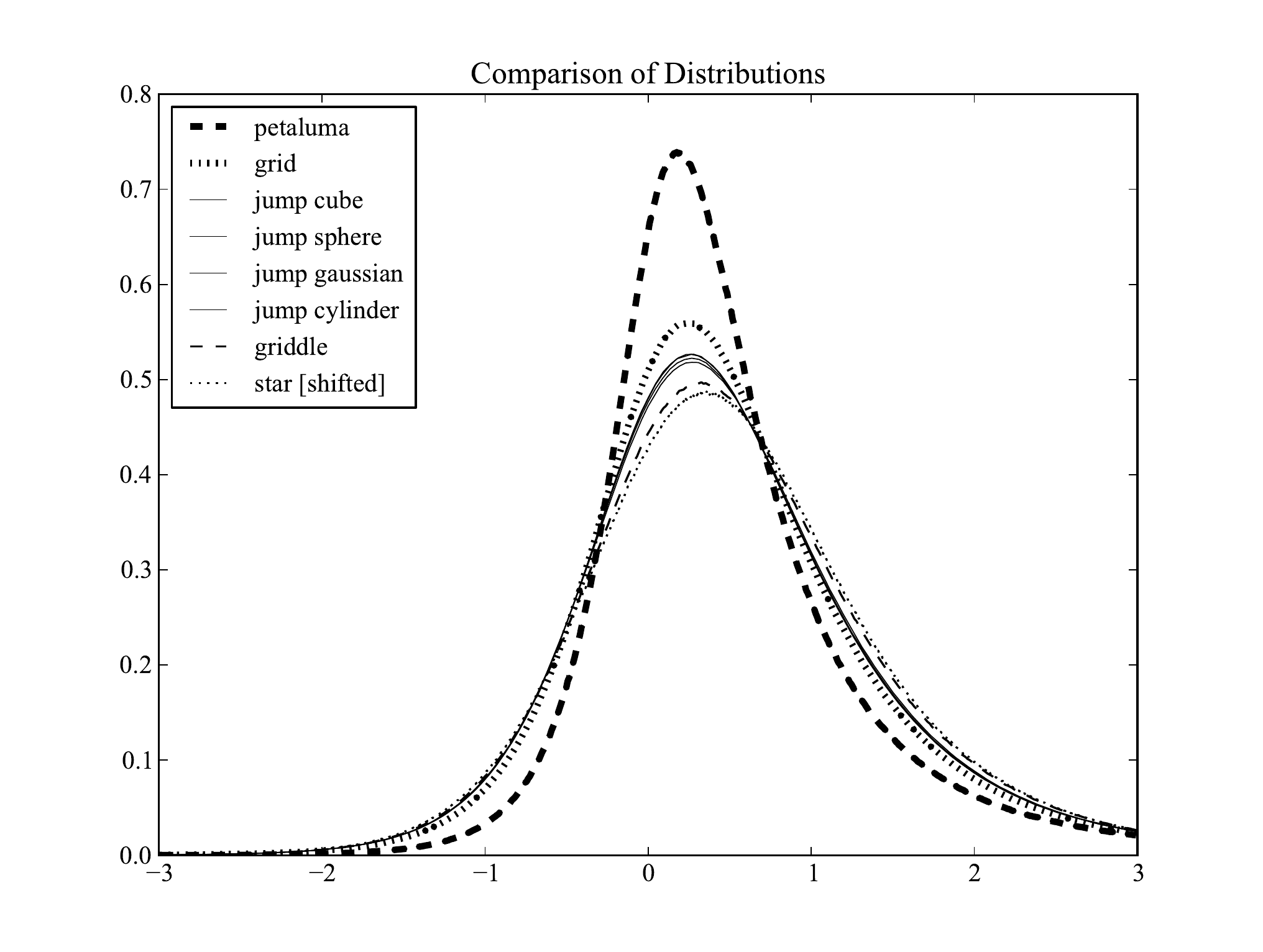}
\captionsetup{width=0.75\textwidth}
\caption{The distribution of $c_2(K_n)/n^2$ in several random knot models, for $n=80$ or $81$, based on $10^8$~random samples each, and normalized to have variance one. Only the star histogram was shifted to compensate for its rightward drift. \\ ~}
\label{compare}
\end{figure}

Our computations and experiments~\cite{even2016invariants, workinprogress} show that the asymptotic distributions of the Casson invariant in models of the third class share several important features. In Figure~\ref{compare}, we exhibit numerically generated histograms of the Casson invariant for three models: Petaluma~(\ref{petaluma-model}), grid~(\ref{grid-model}), and several random jump models~(\ref{jump-model}). They all seem to converge to continuous unimodal limit distributions on $\R$, with two-sided exponentially decaying tails, strictly positive expectations and similarly asymmetric shapes.

Even though models of the second class also seem to converge to distributions of similar shapes around their expectations, their main terms are inconsistent. In the griddle~(\ref{curves-model}) model $E[c_2]/\sqrt{V[c_2]} = \Theta(1)$, while in the star~(\ref{cross-model}) model $E[c_2]/\sqrt{V[c_2]} = \Theta(n)$.

We hope that extending such comparisons to other invariants would shed more light on the above questions of classification and universality.

\bigskip

\newpage

\paragraph{Open Problems}\label{open-section}

Models of the third class outlined above seem especially interesting from a knot-theoretic point of view. They presumably avoid phenomena of local knotting or ``flatness'', and their finite type invariants seem to follow well-behaved distributions. 

We close our review by listing some of the desired features of these random models, which are yet to be established. 

\begin{conj*}
Let $K_n$ be a random knot, sampled from any of the following models: Random Jump~(\ref{jump-model}), Petaluma~(\ref{petaluma-model}), Grid~(\ref{grid-model}). Then,
\begin{itemize}
\setlength\itemsep{0.25em}
\item
With high probability $K_n$ is prime, and even hyperbolic. 
\item
With high probability $K_n$ is non-alternating.
\item
The typical crossing number is super-linear: $E[c(K_n)] = \omega(n)$.
\item
The probability of every knot $K$ is sub-exponential: $P[K_n{=}K] = e^{-\omega(n)}$. 
\item
Any finite type invariant of order $m$ has typical order of magnitude $n^m$.
\end{itemize}
\end{conj*}

\paragraph{Implementation Details}

We include here some information about the numerical results that are firstly reported in this paper.

The generation of random knots in various models was performed by a \emph{C++} program, available at~\cite{even2017ftisampler}. The computation of finite type invariants, as in Sections~\ref{fish-section} and~\ref{casson-section}, was carried out using \emph{Gauss diagram formulas}~\cite{chmutov2012introduction}, which can be evaluated in polynomial time. The computations were distributed on up to 168 processors in the computing facilities of the School of Computer Science and Engineering at HUJI. They were supported by ERC 339096.

The formulas for invariants of random grid and griddle knots with $2n$ segments in Sections~\ref{grid-model} and~\ref{curves-model}, were derived by automated case analysis of the many possible configurations of the involved crossings. It was implemented in a \emph{Python} program, available at~\cite{even2017abcdefg}. These computations took several hours on a~PC.

The data in Figure~\ref{vol} was obtained from the \emph{Sage} software \emph{SnapPy}~\cite{snappy}, that approximates the hyperbolic volume of a link by finding a triangulation of its complement with compatible hyperbolic structure. In order to make the random samples suitable as input for the program, we first converted them from petal diagrams to braids, as shown in Figure~\ref{petalstarbraid}. Some concerns regarding the verification of hyperbolicity and the stability of the computed volume are discussed by Kehne~\cite{uberluma}. Our results are available, together with the source code that generated them, at~\cite{even2017volume}. The computation took several days on a~PC.

\medskip

\section*{Acknowledgment}

I wish to thank Joel Hass, Nati Linial, and Tahl Nowik for many hours of insightful discussions on knots, random knots and many other knotty subjects during our joint work. Their knowledge, vision and ideas are hopefully reflected in this article. All false conjectures are mine.

I also wish to thank Robert Adler, Eric Babson, Rami Band, Itai Benjamini, Harrison Chapman, Moshe Cohen, Nathan Dunfield, Dima Jakobson, Sunder Ram Krishnan, Gal Lavi, Neal Madras, Igor Rivin, Dani Wise, and many others for valuable discussions on randomized knot models. Finally, I thank Moshe Cohen and the reviewers for their helpful suggestions.

\medskip

This article is based on the introduction to the author's doctoral thesis at the Hebrew University of Jerusalem under the supervision of Nati Linial. It was supported by BSF Grant 2012188. 

This final version was completed while the author was a postdoctoral fellow at the Institute for Computational and Experimental Research in Mathematics (ICERM) during the Fall of 2016. The author was supported in part by the Rothschild fellowship.

\bigskip

{\noindent \bfseries The author states that there is no conflict of interest. }

\bigskip

{

\footnotesize
\bibliographystyle{alpha}
\bibliography{models}
}

\end{document}